\documentclass[11pt,twoside]{article}
\usepackage[T1]{fontenc}
\usepackage{amsmath,amsfonts,amssymb,amsthm,fullpage,bbm}
\usepackage{fancyhdr,graphicx,color}
\usepackage{epsfig}
\usepackage{yhmath}

\newtheorem{thm}{Theorem}

\newtheorem{lem}{Lemma}

\theoremstyle{definition}

\theoremstyle{remark}
\newtheorem{rem}{Remark}

\DeclareMathOperator{\cyl}{cyl}

\DeclareMathOperator{\di}{div}
\DeclareMathOperator{\disc}{disc}
\DeclareMathOperator{\hyp}{hyp}
\DeclareMathOperator{\diam}{diam}

\DeclareMathOperator{\tub}{tub}

\DeclareMathOperator{\flow}{flow}
\DeclareMathOperator{\deter}{det}

\newcommand{\eps}{\varepsilon}

\def\PP{\mathbb{P}}
\def\RR{\mathbb{R}}
\def\EE{\mathbb{E}}
\def\NN{\mathbb{N}}
\def\ZZ{\mathbb{Z}}
\def\E{\mathcal{E}}

\def\H{\mathcal{H}}
\def\N{\mathcal{N}}
\def\I{\mathcal{I}}

\def\G{\Gamma}
\def\O{\Omega}
\def\p{\partial}
\def\L{\mathcal{L}}
\def\P{\mathcal{P}}
\def\B{\mathcal{B}}
\def\C{\mathcal{C}}

\def\R{\mathcal{R}}
\def\V{\mathcal{V}}
\def\U{\mathcal{U}}
\def\M{\mathcal{M}}
\def\S{\mathcal{S}}

\def\ga{\gamma}
\def\al{\alpha}
\def\de{\delta}

\def\Pco{\overset{\circ}{\wideparen{\RR^d\smallsetminus P}}}
\def\Po{\overset{\circ}{P}}

\newcommand{\rb}{\partial^{ *}}
\newcommand{\vol}{\mathcal{L}^d}
\newcommand{\ld}{\mathcal{L}^d}
\newcommand{\bGa}{\overline{\Gamma}}
\newcommand{\bro}{\smash{{B}^{\!\!\!\!\raise5pt\hbox{\scriptsize o}}}}
\newcommand{\dro}{\smash{{D}^{\!\!\!\!\raise5pt\hbox{\scriptsize o}}}}
\newcommand{\lro}{\smash{{L}^{\!\!\!\!\raise5pt\hbox{\scriptsize o}}}}
\newcommand{\qro}{\smash{{Q}^{\!\!\!\!\raise5pt\hbox{\scriptsize o}}}}
\newcommand{\uro}{\smash{{U}^{\!\!\!\!\raise5pt\hbox{\scriptsize o}}}}
\newcommand{\ero}{\smash{{E}^{\!\!\!\!\raise5pt\hbox{\scriptsize o}}}}
\newcommand{\cero}{\smash{{\calE}^{\!\!\!\!\raise5pt\hbox{\scriptsize o}}}}
\newcommand{\cgcirc}{\smash{{\calG}^{\!\!\!\raise5pt\hbox{\scriptsize o}}}}
\newcommand{\cdero}{\smash{{\mathbb{E}}^{\!\!\!\!\raise5pt\hbox{\scriptsize o}}}}
\newcommand{\bgcirc}{\smash{{\mathbb{G}}^{\!\!\!\!\raise5pt\hbox{\scriptsize o}}}}
\newcommand{\begcirc}{{{{\mathbb E}\cap\GGG}^{
\!\!\!\!\!\!\!\!\!\raise3pt\hbox{\scriptsize o}}}}

\newcommand{\aro}{\smash{{A}^{\!\!\!\raise5pt\hbox{\scriptsize o}}}}
\newcommand{\thro}{\smash{{\Th}^{\!\!\!\!\!\raise5.2pt\hbox{\scriptsize
        o}}}}

\renewcommand{\qed}{$\ \ \ \Box$}

\title{\huge Law of large numbers for the maximal flow through a domain
  of $\RR^d$ in first passage percolation} 
\author{}
\date{}

\begin{document}
\maketitle

\thispagestyle{empty}

\begin{center}
\vskip-1cm {\Large Rapha\"el Cerf}\\
{\it Universit\'e Paris Sud, Laboratoire de Math\'ematiques, b\^atiment 425\\
91405 Orsay Cedex, France}\\
{\it E-mail:} rcerf@math.u-psud.fr\\
\vskip0.5cm and\\
\vskip0.5cm {\Large Marie Th\'eret}\\
{\it \'Ecole Normale Sup\'erieure, D\'epartement Math\'ematiques et
Applications, 45 rue d'Ulm\\ 75230 Paris Cedex 05, France}\\
{\it E-mail:} marie.theret@ens.fr
\end{center}

\noindent
{\bf Abstract:}
We consider the standard first passage percolation model in the rescaled graph
$\ZZ^d/n$ for $d\geq 2$, and a domain $\O$ of boundary $\G$ in
$\RR^d$. Let $\G^1$ and $\G^2$ be two disjoint open subsets of $\G$, representing
the parts of $\G$ through which some water can enter and escape
from $\O$. We investigate the asymptotic behaviour of the flow $\phi_n$
through a discrete version $\O_n$ of $\O$ between the corresponding discrete
sets $\G^1_n$ and $\G^2_n$. We prove that under some
conditions on the regularity of the domain and on the law of the capacity of
the edges, $\phi_n$ converges almost surely towards a constant
$\phi_{\O}$, which is the solution of a continuous non-random min-cut
problem. Moreover, we give a necessary and sufficient condition on the
law of the capacity of the edges to ensure that $\phi_{\O} >0$.\\




\section{First definitions and main result}

We use many notations introduced in \cite{Kesten:StFlour} and
\cite{Kesten:flows}. Let $d\geq2$. We consider the graph $(\mathbb{Z}^{d}_n,
\mathbb E ^{d}_n)$ having for vertices $\mathbb Z ^{d}_n = \ZZ^d/n$ and for edges
$\mathbb E ^{d}_n$, the set of pairs of nearest neighbours for the standard
$L^{1}$ norm. With each edge $e$ in $\mathbb{E}^{d}_n$ we associate a random
variable $t(e)$ with values in $\mathbb{R}^{+}$. We suppose that the family
$(t(e), e \in \mathbb{E}^{d}_n)$ is independent and identically distributed,
with a common law $\Lambda$: this is the standard model of
first passage percolation on the graph $(\mathbb{Z}^d_n,
\mathbb{E}^d_n)$. We interpret $t(e)$ as the capacity of the edge $e$; it
means that $t(e)$ is the maximal amount of fluid that can go through the
edge $e$ per unit of time.

\begin{figure}[!ht]
\centering
\begin{picture}(0,0)%
\includegraphics{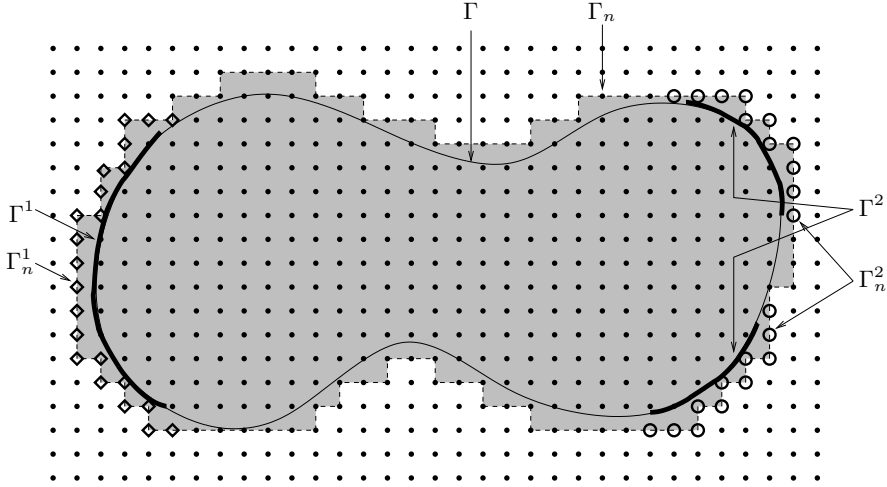}%
\end{picture}%
\setlength{\unitlength}{1973sp}%
\begingroup\makeatletter\ifx\SetFigFont\undefined%
\gdef\SetFigFont#1#2#3#4#5{%
  \reset@font\fontsize{#1}{#2pt}%
  \fontfamily{#3}\fontseries{#4}\fontshape{#5}%
  \selectfont}%
\fi\endgroup%
\begin{picture}(10380,6039)(1261,-6694)
\put(11626,-3361){\makebox(0,0)[lb]{\smash{{\SetFigFont{9}{10.8}{\rmdefault}{\mddefault}{\updefault}{\color[rgb]{0,0,0}$\Gamma^2$}%
}}}}
\put(1276,-3436){\makebox(0,0)[rb]{\smash{{\SetFigFont{9}{10.8}{\rmdefault}{\mddefault}{\updefault}{\color[rgb]{0,0,0}$\Gamma^1$}%
}}}}
\put(1276,-4036){\makebox(0,0)[rb]{\smash{{\SetFigFont{9}{10.8}{\rmdefault}{\mddefault}{\updefault}{\color[rgb]{0,0,0}$\Gamma^1_n$}%
}}}}
\put(11626,-4261){\makebox(0,0)[lb]{\smash{{\SetFigFont{9}{10.8}{\rmdefault}{\mddefault}{\updefault}{\color[rgb]{0,0,0}$\Gamma^2_n$}%
}}}}
\put(6751,-886){\makebox(0,0)[b]{\smash{{\SetFigFont{9}{10.8}{\rmdefault}{\mddefault}{\updefault}{\color[rgb]{0,0,0}$\Gamma$}%
}}}}
\put(8401,-886){\makebox(0,0)[b]{\smash{{\SetFigFont{9}{10.8}{\rmdefault}{\mddefault}{\updefault}{\color[rgb]{0,0,0}$\Gamma_n$}%
}}}}
\end{picture}%
\caption{Domain $\Omega$.}
\label{chapitre7domaine2}
\end{figure}

We consider an open bounded connected subset $\O$ of $\RR^d$ such that
the boundary $\G = \p \O$ of $\O$ is piecewise
of class $\C^1$ (in particular $\G$ has finite area: $\H^{d-1}(\G)
<\infty$). It means that $\G$ is included in the union of a finite number of
hypersurfaces of class $\C^1$, i.e., in the union of a finite number of
$C^1$ submanifolds of~$\RR^d$ of codimension~$1$. Let $\G^1$, $\G^2$ be two disjoint subsets of $\G$ that are
open in $\G$
We want to define the maximal flow from $\G^1$ to $\G^2$ through $\O$ for
the capacities $(t(e), e\in \EE^d_n)$. We consider a discrete version
$(\O_n, \G_n, \G^1_n, \G^2_n)$ of $(\O, \G, \G^1,\G^2)$ defined by:
$$ \left\{ \begin{array}{l} \O_n \,=\, \{ x\in \ZZ^d_n \,|\,
    d_{\infty}(x,\O) <1/n  \}\,,\\ \G_n \,=\, \{ x\in
    \O_n \,|\, \exists y \notin \O_n\,,\,\, \langle x,y \rangle \in \EE^d_n
    \}\,,\\ \G^i_n \,=\, \{ x\in \G_n \,|\, d_\infty (x, \G^i) <1/n \,,\,\,
    d_\infty (x, \G^{3-i}) \geq 1/n \} \textrm{ for }i=1,2 \,,  \end{array}  \right.  $$
where $d_{\infty}$ is the $L^{\infty}$-distance, the notation $\langle
x,y\rangle$ corresponds to the edge of endpoints $x$ and $y$ (see
figure \ref{chapitre7domaine2}).

We shall study the maximal flow from $\G^1_n$ to $\G^2_n$ in $\O_n$.
Let us define
properly  the maximal flow $\phi(F_1 \rightarrow F_2 \textrm{ in } C)$ from
$F_1$ to $F_2$ in $C$, for $C \subset \mathbb{R}^d$ (or by commodity the
corresponding graph $C\cap \mathbb{Z}^d/n$). We will say that an edge
$e=\langle x,y\rangle$ belongs to a subset $A$ of $\mathbb{R}^{d}$, which
we denote by $e\in A$, if the interior of the segment joining $x$ to $y$ is included in $A$. We define
$\widetilde{\mathbb{E}}_n^{d}$ as the set of all the oriented edges, i.e.,
an element $\widetilde{e}$ in $\widetilde{\mathbb{E}}_n^{d}$ is an ordered
pair of vertices which are nearest neighbours. We denote an element $\widetilde{e} \in \widetilde{\mathbb{E}}_n^{d}$ by $\langle \langle x,y \rangle \rangle$, where $x$, $y \in \mathbb{Z}_n^{d}$ are the endpoints of $\widetilde{e}$ and the edge is oriented from $x$ towards $y$. We consider the set $\mathcal{S}$ of all pairs of functions $(g,o)$, with $g:\mathbb{E}_n^{d} \rightarrow \mathbb{R}^{+}$ and $o:\mathbb{E}_n^{d} \rightarrow \widetilde{\mathbb{E}}_n^{d}$ such that $o(\langle x,y\rangle ) \in \{ \langle \langle x,y\rangle \rangle , \langle \langle y,x \rangle \rangle \}$, satisfying:
\begin{itemize}
\item for each edge $e$ in $C$ we have
$$0 \,\leq\, g(e) \,\leq\, t(e) \,,$$
\item for each vertex $v$ in $C \smallsetminus (F_1\cup F_2)$ we have
$$ \sum_{e\in C\,:\, o(e)=\langle\langle v,\cdot \rangle \rangle}
  g(e) \,=\, \sum_{e\in C\,:\, o(e)=\langle\langle \cdot ,v \rangle
    \rangle} g(e) \,, $$
\end{itemize}
where the notation $o(e) = \langle\langle v,. \rangle \rangle$ (respectively $o(e) = \langle\langle .,v \rangle \rangle$) means that there exists $y \in \mathbb{Z}_n^d$ such that $e = \langle v,y \rangle$ and $o(e) = \langle\langle v,y \rangle \rangle$ (respectively $o(e) = \langle\langle y,v \rangle \rangle$).
A couple $(g,o) \in \mathcal{S}$ is a possible stream in $C$ from
$F_1$ to $F_2$: $g(e)$ is the amount of fluid that goes through the edge $e$, and $o(e)$ gives the direction in which the fluid goes through $e$. The two conditions on $(g,o)$ express only the fact that the amount of fluid that can go through an edge is bounded by its capacity, and that there is no loss of fluid in the graph. With each possible stream we associate the corresponding flow
$$ \flow (g,o) \,=\, \sum_{ u \in F_2 \,,\,  v \notin C \,:\, \langle
  u,v\rangle \in \mathbb{E}_n^{d}} g(\langle u,v\rangle) \mathbbm{1}_{o(\langle u,v\rangle) = \langle\langle u,v \rangle\rangle} - g(\langle u,v\rangle) \mathbbm{1}_{o(\langle u,v\rangle) = \langle\langle v,u \rangle\rangle} \,. $$
This is the amount of fluid that crosses $C$ from $F_1$
  to $F_2$ if the fluid respects the stream $(g,o)$. The maximal flow through
  $C$ from $F_1$ to $F_2$ is the supremum of this quantity over all possible choices of streams
$$ \phi (F_1 \rightarrow F_2 \textrm{ in }C)  \,=\, \sup \{ \flow (g,o)\,|\,
  (g,o) \in \mathcal{S} \}  \,.$$

We denote by
$$ \phi_n \,=\, \phi (\G^1_n \rightarrow \G^2_n \textrm{ in } \O_n) $$
the maximal flow from $\G^1_n$ to $\G^2_n$ in $\O_n$. We will investigate
the asymptotic behaviour of $\phi_n/n^{d-1}$ when $n$ goes to infinity. More precisely, we will show that $(\phi_n/n^{d-1})_{n\geq 1}$
converges towards a constant $\phi_{\O}$ (depending on $\O$,
$\G^1$, $\G^2$, $\Lambda$ and
$d$) when $n$ goes to infinity, and that this constant is strictly
positive if and only if $\Lambda (0) < 1-p_c(d)$, where $p_c(d)$ is
the critical parameter for the bond percolation on $\ZZ^d$. The description of $\phi_{\O}$ will be given in section
\ref{chapitre7deflimite}. Here we state the precise theorem:
\begin{thm}
\label{chapitre7lgn}
We suppose that $\O$ is a Lipschitz domain and that $\G$ is included in the union of a
finite number of oriented hypersurfaces $\S_1,...,\S_r$ of class $\C^1$
which are transverse to each other. We also suppose that $\G^1$ and $\G^2$
are open in $\G$, that their relative boundaries $\p_{\G} \G^1$ and
$\p_{\G}\G^2$ in $\G$ have null $\H^{d-1}$ measure, and that $d(\G^1,
\G^2)>0$.
We suppose that the law $\Lambda$ of the capacity of an edge admits an exponential moment:
$$\exists \theta >0 \qquad \int_{\RR^+} e^{\theta x} d\Lambda (x) \,<\, +\infty \,. $$
Then there exists a finite constant $\phi_{\O} \geq 0$ such that
$$\lim_{n\rightarrow \infty} \frac{\phi_n}{n^{d-1}} \,=\, \phi_{\O} \quad
\textrm{ a.s.}  $$
Moreover, this equivalence holds:
$$ \phi_{\O}\,>\, 0 \, \iff \, \Lambda(0)<1-p_c(d)\,.  $$
\end{thm}

\begin{rem}
In the two companion papers \cite{CerfTheret09inf} and \cite{CerfTheret09sup}, we prove in fact that
the lower large deviations of $\phi_n/n^{d-1}$ below $\phi_{\O}$ are of
surface order, and that the upper large deviations of $\phi_n /
n^{d-1}$ above $\phi_{\O}$ are of volume order (see section \ref{secres} where
these results are presented).
\end{rem}


\section{Computation of $\phi_{\O}$}
\label{chapitre7deflimite}

\subsection{Geometric notations}

\begin{figure}[!ht]
\centering
\begin{picture}(0,0)%
\includegraphics{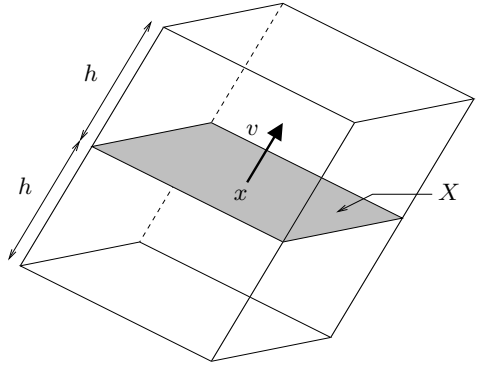}%
\end{picture}%
\setlength{\unitlength}{1973sp}%
\begingroup\makeatletter\ifx\SetFigFont\undefined%
\gdef\SetFigFont#1#2#3#4#5{%
  \reset@font\fontsize{#1}{#2pt}%
  \fontfamily{#3}\fontseries{#4}\fontshape{#5}%
  \selectfont}%
\fi\endgroup%
\begin{picture}(5874,4524)(2539,-6973)
\put(3676,-3436){\makebox(0,0)[rb]{\smash{{\SetFigFont{9}{10.8}{\rmdefault}{\mddefault}{\updefault}{\color[rgb]{0,0,0}$h$}%
}}}}
\put(2851,-4861){\makebox(0,0)[rb]{\smash{{\SetFigFont{9}{10.8}{\rmdefault}{\mddefault}{\updefault}{\color[rgb]{0,0,0}$h$}%
}}}}
\put(5701,-4111){\makebox(0,0)[rb]{\smash{{\SetFigFont{9}{10.8}{\rmdefault}{\mddefault}{\updefault}{\color[rgb]{0,0,0}$v$}%
}}}}
\put(5476,-4936){\makebox(0,0)[b]{\smash{{\SetFigFont{9}{10.8}{\rmdefault}{\mddefault}{\updefault}{\color[rgb]{0,0,0}$x$}%
}}}}
\put(7951,-4936){\makebox(0,0)[lb]{\smash{{\SetFigFont{9}{10.8}{\rmdefault}{\mddefault}{\updefault}{\color[rgb]{0,0,0}$X$}%
}}}}
\end{picture}%
\caption{Cylinder $\cyl(X,h)$.}
\label{chapitre7cylindrerect}
\end{figure}

We start with some geometric definitions. For a subset $X$ of
$\mathbb{R}^d$, we denote by $\mathcal{H}^s (X)$ the $s$-dimensional
Hausdorff measure of $X$ (we will use $s=d-1$ and $s=d-2$). The
$r$-neighbourhood $\V_i(X,r)$ of $X$ for the distance $d_i$, that can be
the Euclidean distance if $i=2$ or the $L^\infty$-distance if $i=\infty$,
is defined by
$$ \V_i (X,r) \,=\, \{ y\in \RR^d\,|\, d_i(y,X)<r\}\,.  $$
If $X$ is a subset of $\RR^d$ included in an hyperplane of $\RR^d$ and of
codimension $1$ (for example a non degenerate hyperrectangle), we denote by
$\hyp(X)$ the hyperplane spanned by $X$, and we denote by $\cyl(X, h)$ the
cylinder of basis $X$ and of height $2h$ defined by
$$ \cyl (X,h) \,=\, \{x+t v \,|\, x\in X \,,\,  t\in
[-h,h]    \}\,,$$
where $v$ is one of the two unit vectors orthogonal to $\hyp(X)$ (see
figure \ref{chapitre7cylindrerect}).
For $x\in
\RR^d$, $r\geq 0$ and a unit vector $v$, we denote by $B(x,r)$ the closed
ball centered at $x$ of radius $r$, by $\disc (x,r,v)$ the
closed disc centered at $x$ of radius $r$ and normal vector $v$, and
by $\hyp (x, v)$ the hyperplane containing $x$ and orthogonal to
$v$. We denote by $\alpha_d$ the volume of a unit ball in $\RR^d$, and
$\alpha_{d-1}$ the $\H^{d-1}$ measure of a unit disc.



\subsection{Flow in a cylinder}

Here are some particular definitions of flows through a box. Let $A$ be a non degenerate hyperrectangle,
i.e., a box of dimension $d-1$ in $\mathbb{R}^d$. All hyperrectangles will be
supposed to be closed in $\mathbb{R}^d$. We denote by
$v$ one of 
the two unit vectors orthogonal to $\hyp (A)$. For $h$ a
positive real number, we consider the cylinder $\cyl(A,h)$.
The set $\cyl(A,h) \smallsetminus \hyp (A)$ has two connected
components, which we denote by $\mathcal{C}_1(A,h)$ and
$\mathcal{C}_2(A,h)$. For $i=1,2$, let $A_i^h$ be
the set of the points in $\mathcal{C}_i(A,h) \cap \mathbb{Z}_n^d$ which have
a nearest neighbour in $\mathbb{Z}_n^d \smallsetminus \cyl(A,h)$:
$$ A_i^h\,=\,\{x\in \mathcal{C}_i(A,h) \cap
\mathbb{Z}_n^d \,|\, \exists y \in \mathbb{Z}_n^d \smallsetminus \cyl(A,h)
\,,\, \langle x,y \rangle \in \EE^d_n\}\,.$$
Let $T(A,h)$ (respectively $B(A,h)$) be the top
(respectively the bottom) of $\cyl(A,h)$, i.e.,
$$ T(A,h) \,=\, \{ x\in \cyl(A,h) \,|\, \exists y\notin \cyl(A,h)\,,\,\,
\langle x,y\rangle \in \mathbb{E}_n^d \textrm{ and }\langle x,y\rangle
\textrm{ intersects } A+hv  \}  $$
and
$$  B(A,h) \,=\, \{ x\in \cyl(A,h) \,|\, \exists y\notin \cyl(A,h)\,,\,\,
\langle x,y\rangle \in \mathbb{E}_n^d \textrm{ and } \langle x,y\rangle
\textrm{ intersects } A-hv  \} \,.$$
For a given realisation $(t(e),e\in \mathbb{E}_n^{d})$ we define the variable
$\tau (A,h) = \tau(\cyl(A,h), v)$ by
$$ \tau(A,h) \,=\,  \tau(\cyl(A,h), v)\,=\, \phi (A_1^h \rightarrow A_2^h
\textrm{ in } \cyl(A,h)) \,,$$
and the variable $\phi(A,h)= \phi(\cyl(A,h), v)$ by
$$ \phi(A,h) \,=\,\phi(\cyl(A,h), v) \,=\, \phi (B(A,h) \rightarrow T(A,h)
\textrm{ in }   \cyl(A,h))\,, $$ 
where $\phi(F_1 \rightarrow F_2 \textrm{ in } C)$ is the maximal
flow from $F_1$ to $F_2$ in $C$, for $C \subset \mathbb{R}^d$ (or by
commodity the corresponding graph $C\cap \mathbb{Z}^d/n$) defined
previously. The dependence in $n$ is implicit here, in fact we can
also write $\tau_n (A,h)$ and $\phi_n(A,h)$ if we want to emphasize
this dependence on the mesh of the graph.


\subsection{Max-flow min-cut theorem}


The maximal flow
$\phi (F_1\rightarrow F_2 \textrm{ in } C)$ can be expressed differently
thanks to the max-flow min-cut theorem (see \cite{Bollobas}). We need some
definitions to state this result.
A path on the graph $\mathbb{Z}_n^{d}$ from $v_{0}$ to $v_{m}$ is a sequence $(v_{0}, e_{1}, v_{1},..., e_{m}, v_{m})$ of vertices $v_{0},..., v_{m}$ alternating with edges $e_{1},..., e_{m}$ such that $v_{i-1}$ and $v_{i}$ are neighbours in the graph, joined by the edge $e_{i}$, for $i$ in $\{1,..., m\}$.
A set $E$ of edges in $C$ is said to cut $F_1$ from $F_2$ in
$C$ if there is no path from $F_1$ to $F_2$ in $C \smallsetminus
E$. We call $E$ an $(F_1,F_2)$-cut if $E$ cuts $F_1$ from $F_2$ in $C$
and if no proper subset of $E$ does. With each set $E$ of edges we
associate its capacity which is the variable
$$ V(E)\, = \, \sum_{e\in E} t(e) \, .$$
The max-flow min-cut theorem states that
$$ \phi(F_1\rightarrow F_2\textrm{ in } C) \, = \, \min \{ \, V(E) \, | \, E
\textrm{ is a } (F_1,F_2)\textrm{-cut} \, \} \, .$$
In fact, as we will see in section \ref{seccontmincut}, $\phi_{\O}$ is
a continuous equivalent of the discrete min-cut.


\subsection{Definition of $\nu$}

The asymptotic behaviour of the rescaled expectation of $\tau_n (A,h)$
for large $n$ is well known, thanks to the almost subadditivity of
this variable. We recall the following result:
\begin{thm}
We suppose that
$$\int_{[0,+\infty[} x \, d\Lambda (x) \,<\,\infty \,.$$
Then for each unit vector $v$ there exists a constant $\nu (d,
\Lambda, v) = \nu(v)$ (the dependence on $d$ and $\Lambda$ is
implicit) such that for every non degenerate hyperrectangle $A$
orthogonal to $v$ and for every strictly positive
constant $h$, we have 
$$ \lim_{n\rightarrow \infty} \frac{\EE [\tau_n(A,h)]}{ n^{d-1}
  \H^{d-1}(A)} \,=\, \nu(v) \,.$$
\end{thm}
For a proof of this proposition, see \cite{RossignolTheret08b}. We
emphasize the fact that the limit depends on the direction of $v$, but
not on $h$ nor on the hyperrectangle $A$ itself.

We recall some geometric properties of the map $\nu: v\in S^{d-1}
\mapsto \nu(v)$, under the only condition on $\Lambda$ that $\EE(t(e))<\infty$. They
have been stated in the section 4.4 of \cite{RossignolTheret08b}. There exists a unit vector $v_0$ such that
$\nu(v_0)=0$ if and only if for all
 unit vector $v$, $\nu(v)=0$, and it happens if and only if $\Lambda(\{0\}) \geq
 1-p_c(d)$. This property has been proved by Zhang in \cite{Zhang}. Moreover, $\nu$ satisfies the weak
 triangle inequality, i.e., if $(ABC)$ is a non degenerate triangle in
 $\mathbb{R}^d$ and $v_A$, $v_B$ and $v_C$ are the
 exterior normal unit vectors to the sides $[BC]$, $[AC]$, $[AB]$ in the
 plane spanned by $A$, $B$, $C$, then
$$ \mathcal{H}^1 ([AB]) \nu(v_C) \,\leq\, \mathcal{H}^1 ([AC])
\nu(v_B) + \mathcal{H}^1 ([BC]) \nu(v_A) \,. $$
This implies that the homogeneous extension $\nu_0$ of $\nu$ to $\RR^d$,
defined by $ \nu_0(0) =0 $ and for all $ w$ in $\RR^d$, 
$$ \nu_0(w)\, =\, |w|_2 \nu (w/|w|_2)\,, $$
is a convex function; in particular, since $\nu_0$ is finite, it is
continuous on $\RR^d$. We denote by $\nu_{\min}$ (respectively
$\nu_{\max}$) the infimum (respectively supremum) of $\nu$ on $S^{d-1}$.


\subsection{Continuous min-cut}
\label{seccontmincut}

\begin{figure}[!ht]
\centering
\begin{picture}(0,0)%
\includegraphics{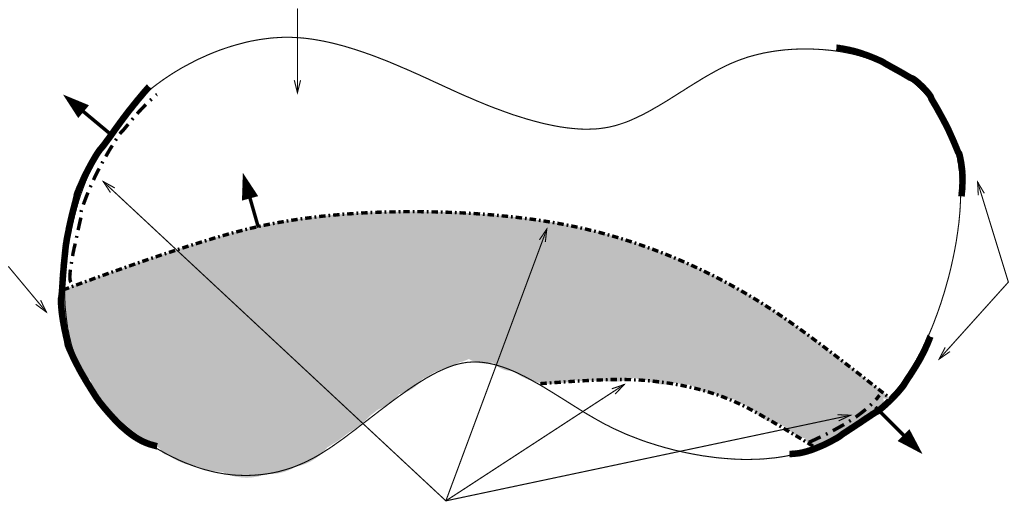}%
\end{picture}%
\setlength{\unitlength}{1973sp}%
\begingroup\makeatletter\ifx\SetFigFont\undefined%
\gdef\SetFigFont#1#2#3#4#5{%
  \reset@font\fontsize{#1}{#2pt}%
  \fontfamily{#3}\fontseries{#4}\fontshape{#5}%
  \selectfont}%
\fi\endgroup%
\begin{picture}(9780,5502)(1411,-6682)
\put(11176,-4261){\makebox(0,0)[lb]{\smash{{\SetFigFont{9}{10.8}{\rmdefault}{\mddefault}{\updefault}{\color[rgb]{0,0,0}$\Gamma^2$}%
}}}}
\put(1426,-4111){\makebox(0,0)[rb]{\smash{{\SetFigFont{9}{10.8}{\rmdefault}{\mddefault}{\updefault}{\color[rgb]{0,0,0}$\Gamma^1$}%
}}}}
\put(4276,-1411){\makebox(0,0)[b]{\smash{{\SetFigFont{9}{10.8}{\rmdefault}{\mddefault}{\updefault}{\color[rgb]{0,0,0}$\Omega$}%
}}}}
\put(3901,-3136){\makebox(0,0)[lb]{\smash{{\SetFigFont{9}{10.8}{\rmdefault}{\mddefault}{\updefault}{\color[rgb]{0,0,0}$v_F(x)$}%
}}}}
\put(3751,-3961){\makebox(0,0)[lb]{\smash{{\SetFigFont{9}{10.8}{\rmdefault}{\mddefault}{\updefault}{\color[rgb]{0,0,0}$x$}%
}}}}
\put(5176,-4411){\makebox(0,0)[lb]{\smash{{\SetFigFont{9}{10.8}{\rmdefault}{\mddefault}{\updefault}{\color[rgb]{0,0,0}$F$}%
}}}}
\put(2026,-2236){\makebox(0,0)[rb]{\smash{{\SetFigFont{9}{10.8}{\rmdefault}{\mddefault}{\updefault}{\color[rgb]{0,0,0}$v_{\Omega}(z)$}%
}}}}
\put(2626,-2986){\makebox(0,0)[lb]{\smash{{\SetFigFont{9}{10.8}{\rmdefault}{\mddefault}{\updefault}{\color[rgb]{0,0,0}$z$}%
}}}}
\put(10351,-5911){\makebox(0,0)[lb]{\smash{{\SetFigFont{9}{10.8}{\rmdefault}{\mddefault}{\updefault}{\color[rgb]{0,0,0}$v_{(F\cap\O)}(y)$}%
}}}}
\put(10201,-5461){\makebox(0,0)[b]{\smash{{\SetFigFont{9}{10.8}{\rmdefault}{\mddefault}{\updefault}{\color[rgb]{0,0,0}$y$}%
}}}}
\put(5776,-6586){\makebox(0,0)[b]{\smash{{\SetFigFont{9}{10.8}{\rmdefault}{\mddefault}{\updefault}{\color[rgb]{0,0,0}$(\partial F \cap \Omega) \cup (\Gamma^2 \cap \partial(F \cap \Omega)) \cup (\Gamma^1 \cap \partial (\Omega \smallsetminus F))$}%
}}}}
\end{picture}%
\caption{The set $(\p F \cap \O) \cup (\G^2 \cap \p (F\cap \O)) \cup (\G^1
  \cap \p(\O\smallsetminus F)) $.}
\label{chapitre7covering2}
\end{figure}

We give here a definition of $\phi_{\Omega}$ and of another constant $\widetilde{\phi_{\O}}$ in
terms of the map $\nu$.
For a subset $F$ of $\RR^d$, we define the perimeter of $F$ in
$\O$ by
$$  \P(F, \O)  \,=\, \sup \left\{\int_F \di f(x) d \L^d(x), \, f\in \C_c^\infty
(\O, B(0,1))  \right\} \,, $$
where $\C_c^\infty (\O, B(0,1))$ is the set of the functions of class
$\C^{\infty}$ from $\RR^d$ to $B(0,1)$, the ball centered at $0$ and of
radius $1$ in $\RR^d$, having a compact support included in $\O$, and $\di$
is the usual divergence operator. The perimeter $\P(F)$ of $F$ is defined as $\P(F, \RR^d)$.
We denote by $\p F$ the boundary of $F$, and by
$\p^* F$ the reduced boundary of $F$. At any point $x$ of $\p^* F$, the set
$F$ admits a unit exterior normal vector $v_F(x)$ at $x$ in a measure
theoretic sense (for
definitions see for
example \cite{Cerf:StFlour} section 13). For all $F \subset \RR^d$ of
finite perimeter in $\O$, we define
\begin{align*}
\I_{\O}(F) &\,=\, \int_{\p^* F \cap \O} \nu(v_F(x)) d\H^{d-1}(x) +
\int_{\G^2 \cap \p^* (F\cap \O)} \nu(v_{(F\cap \O)}(x)) d\H^{d-1}(x)\\ & \qquad \qquad +
\int_{\G^1 \cap \p^* (\O\smallsetminus F)} \nu(v_{\O}(x)) d\H^{d-1}(x)\,.
\end{align*}
If $\P(F, \O) = +\infty$, we define $\I_{\O}(F) =+\infty$. Finally, we define
$$ \phi_{\O} \,=\, \inf \{ \I_{\O}(F)\,|\, F\subset \RR^d  \} 
 \,=\, \inf \{ \I_{\O}(F)\,|\, F\subset \O  \}\,. $$
In the case where $\partial F$ is $\C^1$, $\I_{\O}(F)$ has the simpler
following expression:
\begin{align*}
\I_{\O}(F) &\,=\, \int_{\p F \cap \O} \nu(v_F(x)) d\H^{d-1}(x) +
\int_{\G^2 \cap \p (F\cap \O)} \nu(v_{(F\cap \O)}(x)) d\H^{d-1}(x)\\ & \qquad \qquad +
\int_{\G^1 \cap \p(\O\smallsetminus F)} \nu(v_{\O}(x)) d\H^{d-1}(x)\,.
\end{align*}
The localization of the set along which the previous integrals are done is
illustrated in figure \ref{chapitre7covering2}.

\begin{figure}[!ht]
\centering
\begin{picture}(0,0)%
\includegraphics{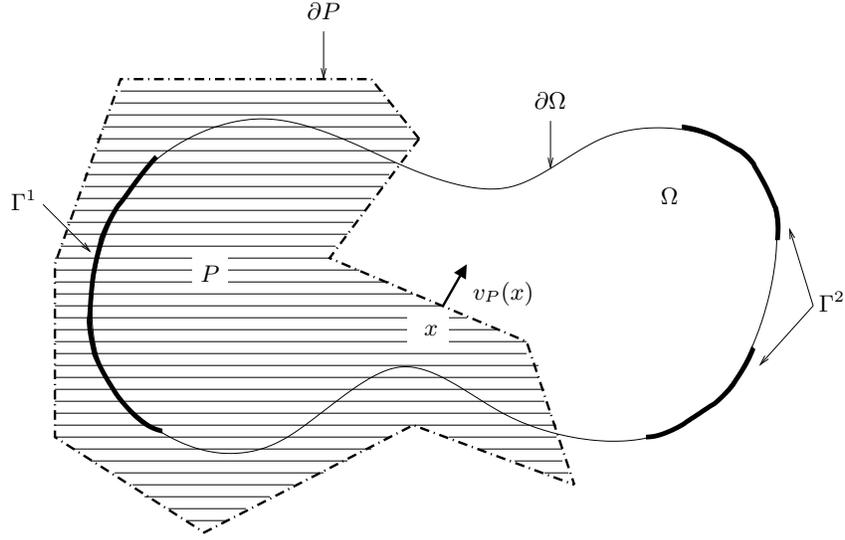}%
\end{picture}%
\setlength{\unitlength}{1973sp}%
\begingroup\makeatletter\ifx\SetFigFont\undefined%
\gdef\SetFigFont#1#2#3#4#5{%
  \reset@font\fontsize{#1}{#2pt}%
  \fontfamily{#3}\fontseries{#4}\fontshape{#5}%
  \selectfont}%
\fi\endgroup%
\begin{picture}(9855,6714)(1336,-7069)
\put(11176,-4261){\makebox(0,0)[lb]{\smash{{\SetFigFont{9}{10.8}{\rmdefault}{\mddefault}{\updefault}{\color[rgb]{0,0,0}$\Gamma^2$}%
}}}}
\put(6826,-4111){\makebox(0,0)[lb]{\smash{{\SetFigFont{9}{10.8}{\rmdefault}{\mddefault}{\updefault}{\color[rgb]{0,0,0}$v_{P}(x)$}%
}}}}
\put(1351,-2986){\makebox(0,0)[rb]{\smash{{\SetFigFont{9}{10.8}{\rmdefault}{\mddefault}{\updefault}{\color[rgb]{0,0,0}$\Gamma^1$}%
}}}}
\put(9301,-2911){\makebox(0,0)[b]{\smash{{\SetFigFont{9}{10.8}{\rmdefault}{\mddefault}{\updefault}{\color[rgb]{0,0,0}$\Omega$}%
}}}}
\put(4951,-586){\makebox(0,0)[b]{\smash{{\SetFigFont{9}{10.8}{\rmdefault}{\mddefault}{\updefault}{\color[rgb]{0,0,0}$\partial P$}%
}}}}
\put(7801,-1711){\makebox(0,0)[b]{\smash{{\SetFigFont{9}{10.8}{\rmdefault}{\mddefault}{\updefault}{\color[rgb]{0,0,0}$\partial \O$}%
}}}}
\put(3526,-3886){\makebox(0,0)[b]{\smash{{\SetFigFont{9}{10.8}{\rmdefault}{\mddefault}{\updefault}{\color[rgb]{0,0,0}$P$}%
}}}}
\put(6301,-4561){\makebox(0,0)[b]{\smash{{\SetFigFont{9}{10.8}{\rmdefault}{\mddefault}{\updefault}{\color[rgb]{0,0,0}$x$}%
}}}}
\end{picture}%
\caption{A polyhedral set $P$ as in the definition of $\widetilde{\phi_{\O}}$.}
\label{chapitre7separantsup}
\end{figure}

When a hypersurface $\S$ is piecewise of
class $\C^1$, we say that $\S$ is
transverse to $\G$ if for all $x\in \S \cap \G$, the
normal unit vectors to $\S$ and $\G$ at $x$ are not collinear; if the
normal vector to $\S$ (respectively to $\G$) at $x$ is not well defined, this property must
be satisfied by all the vectors which are limits of normal unit vectors to
$\S$ (respectively $\G$) at $y \in \S$ (respectively $y\in \G$) when we
send $y$ to $x$ - there is at most a finite number of such limits. We say
that a subset $P$
of $\RR^d$ is polyhedral if its boundary $\p P $ is included in the
union of a finite number of hyperplanes. For each point $x$ of
such a set $P$ which is on the interior of one face of $\p P$, we denote by
$v_{P}(x)$ the exterior unit vector orthogonal to $P$ at $x$. For $A
\subset \RR^d$, we denote by $\overset{\circ}{A}$ the interior of $A$. We define
$\widetilde{\phi_{\O}}$ by
$$ \widetilde{\phi_{\O}} \,=\, \inf \left\{ \I_{\O} (P) \,\Bigg|
  \, \begin{array}{c} P \subset \RR^d\,,\, \overline{\G^1} \subset
\Po \,,\, \overline{\G^2} \subset \Pco\\
 P\textrm{ is polyhedral} \,,\, \p P \textrm{ is transverse to
 }\G \end{array} \right\} \,.$$
Notice that if $P$ is a set such that
$$  \overline{\G^1} \,\subset\, \Po \quad \textrm{and}\quad \overline{\G^2} \subset \Pco \,,$$
then
$$ \I_{\O} (P) \,=\, \int_{\p  P \cap \O} \nu(v_{P}(x)) d \H^{d-1} (x)\,.$$
See figure \ref{chapitre7separantsup} to have an example of such a
polyhedral set $P$.

The definitions of the constants $\phi_{\O}$ and $\widetilde{\phi_{\O}}$
are not very intuitive. We propose to define the notion of a continuous
cutset to have a better understanding of these constants. We say that $\S
\subset \RR^d$ cuts $\G^1$ from $\G^2$ in $\overline{\O}$ if
every continuous path from $\G^1$ to $\G^2$ in $\overline{\O}$ intersects
$\S$. In fact, if $P$ is a polyhedral set of $\RR^d$ such that 
$$\overline{\G^1} \,\subset\, \Po \quad \textrm{and} \quad
\overline{\G^2}\, \subset\, \Pco \,, $$
then $\p P \cap \overline{\O}$ is a continuous cutset from $\G^1$ to $\G^2$
in $\overline{\O}$.
Since $\nu(v)$ is the average
amount of fluid that can cross a
hypersurface of area one in the direction $v$ per
unit of time, it can be interpreted as the capacity of a unitary
hypersurface. Thus $\I_{\O}(P)$ can be interpreted as the capacity of the
continuous cutset $\partial P \cap \overline{\O}$. The constant
$\widetilde{\phi_{\O}}$ is the solution of a min cut problem, because it is equal to the infimum of the capacity of a continuous cutset that
satisfies some specific properties. We can define two other
constants, that are solutions of possibly more intuitive min cuts problems. If
$\S$ is a hypersurface which is piecewise of class $\C^1$, we denote by
$v_{\S}(x)$ one of the two normal unit vectors to $\S$ at $x$ for every
point $x$ at which $\S$ is regular. The $\H^{d-1}$ measure of the points at
which $\S$ is not regular is null. We define
$$\widehat{\phi_{\O}} \,=\, \inf \left\{ \int_{\S \cap \overline{\O}} \nu(v_{S}(x))d\H^{d-1}(x) \,\Bigg|
  \, \begin{array}{c} \S \textrm{ hypersurface piecewise of class }\C^1\\
    \S \textrm{ cuts }\G^1\textrm{ from }\G^2\textrm{ in }\overline{\O} \end{array} \right\}   $$
and
$$\widetriangle{\phi_{\O}} \,=\, \inf \left\{ \int_{\S \cap \overline{\O}} \nu(v_{S}(x))d\H^{d-1}(x) \,\Bigg|
  \, \begin{array}{c} \S \textrm{ polyhedral hypersurface} \\ \S \textrm{
      cuts }\G^1\textrm{ from }\G^2\textrm{ in }\overline{\O} \end{array}
\right\} \,.   $$
We remark that by definition,
$$ \widehat{\phi_{\O}} \,\leq\, \widetriangle{\phi_{\O}} \,\leq\,
\widetilde{\phi_{\O}}\,. $$
We claim that $\phi_{\O} \leq \widehat{\phi_{\O}}$. Let $\S$ be a hypersurface
which is piecewise of class $\C^1$, which cuts $\G^1$ from $\G^2$ in
$\overline{\O}$, and such that 
$$ \int_{\S \cap \overline{\O}} \nu(v_{\S}(x))d\H^{d-1}(x) \,\leq\,  \widehat{\phi_{\O}} + \eta $$
for some positive $\eta$. Let $F$ be the set of the points of
$\overline{\O}\smallsetminus \S$ that can be joined to a point of $\G^1$ by
a continuous path. Then
$$ (\p F \cap \O) \cup (\G^1 \cap \p (\O \smallsetminus F)) \cup (\G^2
\cap \p (F\cap \O)) \,\subset\, \S \cap \overline{\O} \,. $$ 
Thus $F$ is of finite
perimeter in $\O$, and $\I_{\O}(F)$ satisfies
$$ \I_{\O}(F) \,\leq\,  \int_{\S \cap \overline{\O}}
\nu(v_{\S}(x))d\H^{d-1}(x) \,\leq\, \widehat{\phi_{\O}} + \eta\,. $$
Thus we have proved that
$$ \phi_{\O} \,\leq\, \widehat{\phi_{\O}} \,\leq\, \widetriangle{\phi_{\O}} \,\leq\,
\widetilde{\phi_{\O}}\,. $$


\section{State of the art}

\subsection{Existing laws of large numbers}

Only in this section, we consider the
standard first passage percolation model on the graph $(\ZZ^d, \EE^d)$
instead of the rescaled graph $(\ZZ^d_n, \EE^d_n)$. We present here some
laws of large numbers that have been proved about maximal flows. 

Using a subadditive argument and concentration inequalities, Rossignol
and Th\'eret have proved in \cite{RossignolTheret08b} that $\tau(nA,
h(n))$ satisfies a law of large numbers:
\begin{thm}[Rossignol and Th\'eret]
We suppose that
$$ \int_{[0,\infty[} x \, d\Lambda (x) \,<\, \infty\, .  $$
For every unit vector $v$, for every non degenerate hyperrectangle $A$ orthogonal to $v$, for
every height function $h: \NN \rightarrow \RR^+$ satisfying
$\lim_{n\rightarrow \infty} h(n) = +\infty$, we have
$$ \lim_{n\rightarrow \infty} \frac{\tau(nA, h(n))}{\H^{d-1}(nA)}
\,=\, \nu(v) \qquad \textrm{in } L^1 \,. $$
Moreover, if the origin of the graph belongs to $A$, or if
$$ \int_{[0,\infty[} x^{1+\frac{1}{d-1}} \, d\Lambda (x) \,<\, \infty \,, $$
then
$$ \lim_{n\rightarrow \infty} \frac{\tau(nA, h(n))}{\H^{d-1}(nA)}
\,=\, \nu(v) \qquad \textrm{a.s.} $$
\end{thm}

Kesten, Zhang, Rossignol and Th\'eret have studied the maximal flow between the top and the
bottom of straight cylinders. Let us denote by $D({\bf k}, m)$ the cylinder
$$ D({\bf k}, m) \,=\, \prod_{i=1}^{d-1} [0,k_i] \times [0,m]\,, $$
where ${\bf k} = (k_1,...,k_{d-1}) \in \RR^{d-1}$. We denote by $\phi({\bf
  k},m)$ the maximal flow in $ D({\bf k}, m)$ from its top $
\prod_{i=1}^{d-1} [0,k_i] \times \{ m \}$ to its bottom $
\prod_{i=1}^{d-1} [0,k_i] \times \{0\}$. Kesten proved in
\cite{Kesten:flows} the following result:
\begin{thm}[Kesten]
Let $d=3$. We suppose that $\Lambda(0)<p_0$ for some fixed $p_0\geq 1/27$, and that
$$ \exists \gamma >0 \qquad \int_{[0,+\infty[} e^{\gamma x} \,d \Lambda(x)
\,<\,\infty\,.  $$
If $m=m({\bf k})$ goes to infinity with $k_1 \geq k_2$ in such a way that
$$ \exists \delta >0 \qquad \lim_{k_1\geq k_2 \rightarrow \infty}
k^{-1+\delta} \log m({\bf k}) \,=\,0\,,  $$
then
$$  \lim_{k_1\geq k_2 \rightarrow \infty} \frac{\phi({\bf
  k},m)}{k_1 k_2} \,=\, \nu((0,0,1)) \qquad \textrm{a.s. and in }L^1\,. $$
Moreover, if $\Lambda(0)>1-p_c(d)$, where $p_c(d)$ is the critical parameter for the
standard bond percolation model on $\ZZ^d$, and if
$$  \int_{[0,+\infty[} x^6 \, d\Lambda(x) \,<\,\infty\,, $$
there exists a constant $C=C(F) <\infty$ such that for all $m=m({\bf k})$
that goes to infinity with $k_1\geq k_2$ and satisfies
$$ \liminf_{k_1\geq k_2 \rightarrow \infty} \frac{m({\bf k})}{k_1 k_2} \,>\,C\,,  $$
 for all $k_1\geq k_2$ sufficiently large, we have
$$ \phi({\bf k},m) \,=\, 0 \qquad \textrm{a.s.} $$
\end{thm}
Zhang improved this result in \cite{Zhang07} where he proved the following
theorem:
\begin{thm}[Zhang]
Let $d\geq 2$. We suppose that
$$ \exists \gamma >0 \qquad \int_{[0,+\infty[} e^{\gamma x} \,d\Lambda (x)\,<\,\infty\,. $$
Then for all $m=m({\bf k})$ that goes to infinity when all the $k_i$,
$i=1,...,d-1$ go to infinity in such a way that
$$ \exists \delta \in ]0,1] \qquad \log m({\bf k})
\,\leq\,\max_{i=1,...,d-1} k_i^{1-\delta}\,, $$
we have
$$ \lim_{k_1,...,k_{d-1} \rightarrow \infty} \frac{\phi({\bf k},
  m)}{\prod_{i=1}^{d-1} k_i} \,=\, \nu((0,...,0,1)) \qquad \textrm{a.s. and
in } L^1\,.  $$
Moreover, this limit is positive if and only if $\Lambda(0)<1-p_c(d)$.
\end{thm}
To show this theorem, Zhang obtains first an important control on the
number of edges in a minimal cutset. Finally, Rossignol and Th\'eret
improved Zhang's result in \cite{RossignolTheret08b} in the particular
case where the dimensions of the basis of the straight cylinder go to
infinity all at the same speed. They obtain the following result:
\begin{thm}[Rossignol and Th\'eret]
We suppose that
$$ \int_{[0,\infty[} x \, d\Lambda (x)\,<\, \infty \,. $$
For every straight hyperrectangle $A = \prod_{i=1}^{d-1} [0,a_i] \times
\{0 \}$ with $a_i >0$ for all $i$, for every height function $h:\NN \rightarrow
\RR^+$ satisfying $\lim_{n\rightarrow \infty} h(n) = +\infty$ and
$\lim_{n\rightarrow \infty} \log h(n) / n^{d-1} =0$, we have
$$ \lim_{n\rightarrow \infty} \frac{\phi(nA,h(n))}{\H^{d-1} (nA)}
\,=\, \nu ((0,...,0,1)) \qquad \textrm{a.s. and in }L^1 \,.  $$
\end{thm}

In dimension two, more results are known. We present here two of
them. Rossignol and Th\'eret have studied in \cite{RossignolTheret09} the maximal flow from the
top to the bottom of a tilted cylinder in dimension two, and they
have proved the following theorem (Corollary 2.10 in \cite{RossignolTheret09}):
\begin{thm}[Rossignol and Th\'eret]
Let $v$ be a unit vector, $A$ a non degenerate line-segment orthogonal to
$v$, $h:\NN \rightarrow \RR^+$ a height function satisfying $\lim_{n
  \rightarrow \infty} h(n) = +\infty$ and
    $\lim_{n\rightarrow \infty} \log h(n) / n =0$. We suppose that there
    exists $\alpha \in [0,\pi/2]$ such that
$$ \lim_{n\rightarrow \infty} \frac{2h(n)}{\H^1(nA)} \,=\, \tan \alpha \,. $$
Then, if 
$$  \int_{[0,\infty[} x \, d\Lambda (x) \,<\, \infty \,,  $$
we have
$$ \lim_{n \rightarrow \infty} \frac{\phi(nA,h(n))}{\H^1(nA)} \,=\,
    \inf \left\{ \frac{\nu(v')}{v\cdot v'} \,\Big\vert\, v' \textrm{ satisfies
      } v\cdot v' \geq \cos \alpha \right\} \qquad \textrm{in }L^1\,. $$
Moreover, if the origin of the graph is the middle of $A$, or if
$$  \int_{[0,\infty[} x^2 \, d\Lambda (x) \,<\, \infty \,,  $$
then we have
$$ \lim_{n \rightarrow \infty} \frac{\phi(nA,h(n))}{\H^1(nA)} \,=\,
    \inf \left\{ \frac{\nu(v')}{v\cdot v'} \,\Big\vert\, v' \textrm{ satisfies
      } v\cdot v' \geq \cos \alpha \right\} \qquad \textrm{a.s.} $$
\end{thm}

Garet studied in \cite{Garet} the maximal flow $\sigma (A)$ between a
convex bounded set $A$ and infinity in the case $d=2$. By an extension of
the max flow - min cut theorem to non finite graphs, Garet proves in
\cite{Garet} that this maximal flow is equal to the minimal capacity of a
set of edges that cuts all paths from $A$ to infinity. Let $\p A$ be the
boundary of $A$, and $\p ^* A$ the set of the points $x \in \p A$ at which
$A$ admits a unique exterior normal unit vector $v_A(x)$ in a measure
theoretic sense (see \cite{Cerf:StFlour}, section 13, for a precise definition). If $A$ is a convex set, the set $\p^* A$ is also equal to
the set of the points $x\in \p A$ at which $A$ admits a unique exterior
normal vector in the classical sense, and this vector is $v_A(x)$. Garet proved
the following theorem:
\begin{thm}[Garet]
Let $d=2$. We suppose that $\Lambda (0) <1-p_c(2) = 1/2$ and that
$$ \exists \gamma >0 \qquad \int_{[0,+\infty[} e^{\gamma x}\, d\Lambda
(x)\,<\,\infty\,.  $$
Then for all convex bounded set $A$ containing $0$ in its interior, we have
$$ \lim_{n\rightarrow \infty} \frac{\sigma(nA)}{n} \,=\, \int_{\p^* A}
\nu(v_A(x)) d\H^1 (x) \,=\, \I(A) \,>\,0 \qquad\textrm{a.s.}  $$
Moreover, for all $\eps >0$, there exist constants $C_1$, $C_2>0$ depending
on $\eps$ and $\Lambda$ such that
$$ \forall n\geq 0 \qquad \PP \left[\frac{\sigma(nA)}{n\I(A)} \notin
  ]1-\eps, 1+\eps[  \right] \,\leq\, C_1 \exp (-C_2 n)\,. $$ 
\end{thm}

Nevertheless, a law of large numbers for the maximal flow from the top
to the bottom of a tilted cylinder for $d\geq 3$ was not proved yet. In fact, the lack of symmetry of
the graph induced by the slope of the box is a major issue to extend
the existing results
concerning straight cylinders to tilted cylinders. The theorem of Garet was
not extended to dimension $d\geq 3$ either. Theorem \ref{chapitre7lgn}
applies to the maximal flow from the top to the bottom of a tilted
cylinder. Thus it is a generalisation of the laws of large numbers of
Kesten, Zhang, Rossignol and Th\'eret for
the variable $\phi$ in straight cylinders, in the particular case
where all the dimensions of the cylinder go to infinity at the same
speed (or, equivalently, the cylinder is fixed and the mesh of the
graph go to zero isotropically). Moreover, it gives a hint of what
could be a generalisation of the result of Garet in higher dimension,
all the more since the expression of the constant $\phi_{\O}$ is a
reminiscent of the value of the limit in Garet's Theorem: the capacity $\I_{\O}$ of a continuous cutset is exactly the
same as the one defined by Garet in \cite{Garet} in dimension two, except
that we consider a maximal flow through a bounded domain, so our capacity
is adapted to deal with specific boundary conditions.

From now on, we work in the rescaled graph $(\ZZ^d_n, \EE^d_n)$.


\subsection{Large deviations for $\phi_n$}
\label{secres}

We present here the two existing results concerning $\phi_n$. We
consider an open bounded connected subset $\O$ of $\RR^d$,
whose boundary $\G$ is piecewise of class $\C^1$, and two disjoint
open subsets $\G^1$ and $\G^2$ of $\G$. The first result 
states that the lower large deviations below $\phi_{\O}$ are of
surface order, and is proved by the authors in \cite{CerfTheret09inf}:
\begin{thm}
\label{chapitre7devinf}
If the law $\Lambda$ of the capacity of an edge admits an exponential moment:
$$\exists \theta >0 \qquad \int_{\RR^+} e^{\theta x} d\Lambda (x) \,<\, +\infty \,, $$
and if $\Lambda (0) < 1-p_c(d)$, then for all $\lambda < \phi_{\O}$,
$$ \limsup_{n\rightarrow \infty} \frac{1}{n^{d-1}}\log \PP [ \phi_n \leq
\lambda n^{d-1}  ] \,<\,0\,.$$
\end{thm}
The second result states that the upper large deviations of $\phi_n$ above
$\widetilde{\phi_{\O}}$ are of volume order and is proved by the
authors in \cite{CerfTheret09sup}:
\begin{thm}
\label{chapitre7devsup}
We suppose that  $d(\G^1, \G^2)>0$.
If the law $\Lambda$ of the capacity of an edge admits an exponential moment:
$$\exists \theta >0 \qquad \int_{\RR^+} e^{\theta x} d\Lambda (x) \,<\, +\infty \,, $$
then for all $\lambda >\widetilde{\phi_{\O}}$,
$$\limsup_{n\rightarrow \infty} \frac{1}{n^d} \log \PP [ \phi_n \geq \lambda n^{d-1}  ] \,<\, 0 \,.$$
\end{thm}

By a simple Borel-Cantelli lemma, these results imply that if $\Lambda$ admits
an exponential moment and if $d(\G^1, \G^2)>0$, then
$$ \phi_{\O} \,\leq\, \liminf_{n\rightarrow \infty} \frac{\phi_n}{n}
\,\leq\, \limsup_{n\rightarrow \infty} \frac{\phi_n}{n}
\,\leq\, \widetilde{\phi_{\O}} \,. $$
Notice here that Theorem \ref{chapitre7devinf} allows us to obtain the
first inequality only under the additional hypothesis that $\Lambda
(0) < 1-p_c(d)$, however if $\Lambda (0) \geq 1-p_c(d)$ we know that $\nu
(v) =0$ for all $v$, so $\phi_{\O}=0$ and the first inequality remains
valid.

Thus, to prove Theorem \ref{chapitre7lgn}, it remains to prove that
$\phi_{\O} = \widetilde{\phi_{\O}}$, and to study
the positivity of $\phi_{\O}$. The equality $\phi_{\O}= \widetilde{\phi_{\O}}$ is a consequence of a
polyhedral approximation of sets having finite perimeter that will be done
in section \ref{chapitre7secap}. The positivity of $\phi_{\O}$ is proved in section \ref{chapitre7secpositif}, using tools
of differential geometry like tubular neighbourhood of paths. These two results are proved by purely geometrical studies. Since the probabilistic part of the proof of Theorem
\ref{chapitre7lgn} is contained in Theorems \ref{chapitre7devinf} and \ref{chapitre7devsup}, we
propose a sketch of the proofs of these two theorems in sections
\ref{seclowersketch} and \ref{secuppersketch} to help the understanding of
the law of large numbers proved in this paper.

Before these two sketches of proofs, we would like to make two
remarks. The first one is that the large deviations that are obtained
in Theorem \ref{chapitre7devinf} and \ref{chapitre7devsup} are of the
relevant order. Indeed, if all the edges in $\O_n$ have a capacity which is abnormally big, then the
maximal flow $\phi_n$ will be abnormally big too. The probability for these
edges to have an abnormally large capacity is of order $\exp -C n^d$ for a
constant $C$, because the number of edges in $\O_n$ is $C' n^{d}$ for a
constant $C'$. On the opposite, if all the edges in a flat layer that
separates $\G^1_n$ from $\G^2_n$ in $\O_n$ have abnormally small capacity,
then $\phi_n$ will be abnormally small. Since the cardinality of such a set
of edges is $D' n^{d-1}$ for a constant $D'$, the probability of this event
is of order $\exp -D n^{d-1}$ for a constant $D$.

The second remark we would like to do is that the condition $d (\G^1,
\G^2) >0$ is relevant in Theorem ~\ref{chapitre7devsup}. First, without
this condition, we cannot be sure that there exists a polyhedral set
$P$ as in the definition of $\widetilde{\phi_{\O}}$, and thus the
polyhedral approximation (see section \ref{chapitre7secap}) cannot be
performed. Moreover, if $d(\G^1, \G^2) =0$, there exists a set of edges
of constant cardinality (not depending on $n$) that 
contains paths from $\G^1_n$ to $\G^2_n$ through $\O_n$ for all $n$ along the
common boundary of $\G^1$ and $\G^2$, and so it may be sufficient for
these edges to have a huge capacity to obtain that $\phi_n$ is abnormally big
too. Thus, we cannot hope to obtain upper large deviations of volume
order (see \cite{Theret:uppertau} for a counter-example). However, we do
not know if this condition is essential for Theorem \ref{chapitre7lgn} to hold.


\subsubsection{Lower large deviations}
\label{seclowersketch}

To prove Theorem \ref{chapitre7devinf}, we have to study the
probability
\begin{equation}
\label{proba1}
 \PP \left[ \phi_n \leq (\phi_{\O} - \eps) n^{d-1}\right]
\end{equation}
for a positive $\eps$. The proof is divided in three steps.\\

{\bf First step:} We consider a set of edges $\E_n$ that cuts
$\Gamma^1_n$ from $\Gamma^2_n$ in $\O_n$, of minimal capacity (so
$\phi_n = V(\E_n)$) and having the minimal number of edges among those
cutsets. We see it as the (edge) boundary of a set $E_n$ which is
included in $\O$. Zhang's estimate of the number of edges in a
minimal cutset (Theorem 1 in \cite{Zhang07}) states that with high
probability, the perimeter $\P (E_n , \O)$ of $E_n$ in $\O$ is
smaller than a constant $\beta$. Thus, $E_n$ belongs to the set
$$ \C_\beta \,=\, \{ F \subset \O \,\vert \, F\subset \O \,,\,\, \P (F,
\O) \leq \beta  \}\,. $$
We endow $\C_\beta$ with the topology $L^1$ associated to the
following distance $d$:
$$ d (F_1,F_2) \,=\, \L^d (F_1 \triangle F_2)\,,$$
where $\L^d$ is the $d$-dimensional Lebesgue measure. For this topology, the set $\C_\beta$ is compact. Thus, if we
associate to each set $F$ in $\C_\beta$ a positive constant $\eps_F$,
and if we denote by $\V(F,\eps_F)$ the neighbourhood of $F$ of radius
$\eps_F$ for the distance $d$ defined above, the collection of these
neighbourhoods is an open covering of $\C_\beta$, and thus by compactness of
$\C_\beta$ we can extract a finite covering:
$$ \exists F_1,...,F_N \qquad \C_\beta \subset \bigcup_{i=1}^N \V(F_i,
\eps_{F_i}) \,.$$
If we find an upper bound on the following probability:
\begin{equation}
\label{proba2}
 \PP \left[ \phi_n \leq (\phi_{\O} - \eps) n^{d-1} \textrm{ and }
  d(E_n,F)\leq \eps_F \right] 
\end{equation}
for each $F$ in $\C_\beta$ and a corresponding $\eps_F$, then we will
obtain an upper bound on the probability (\ref{proba1}).\\

{\bf Second step:} We consider a fixed set $F$ in $\C_\beta$, and we
want to evaluate the probability (\ref{proba2}). So we suppose that
$E_n$ is close to $F$ for the distance $d$, we denote it by $E_n
\approx F$ to simplify the notations. We skip here all the problems of
boundary conditions that arise in the proof of Theorem~\ref{chapitre7devinf}: we suppose that $\I_{\O} (F)$ is equal to the
integral of $\nu$ along $\partial ^* F \cap \O$.

We make a zoom along $\partial F$. Using the Vitali covering Theorem
(Theorem \ref{chapitre7vitali} in section \ref{chapitre7secap}),
we know that there exists a finite number of disjoint balls
$B_j=B(x_j,r_j)$ for $j=1,...,\N$ with $x_j \in \partial F$ such that $\partial F$ is ''almost
flat'' in each ball, and the part of $\partial F$ that is missing in
the covering has a very small area. We denote by $v_j$ the exterior
normal unit vector of $F$ at $x_j$ (we suppose that it exists). Here
''almost flat'' means that
\begin{itemize}
\item[(i)] the capacity of $\partial F$ inside $B_j$ is very close to the
capacity of the flat disc $\hyp (x_j, v_j)  \cap B_j$, i.e., very
close to $\alpha_{d-1} r_j^{d-1} \nu(v_j)$ ;
\item[(ii)] $F\cap B_j \approx B_j^-$, where $B_j^-$ is the lower half part of the
ball $B_j$ in the direction given by $v_j$:
$$ B_j^- \,=\, \{ y\in B_j \,\vert\, (y-x_j) \cdot v_j <0  \}\,. $$
\end{itemize}
Thanks to property (i) and the fact that only a very small area of
$\partial F$ is missing in the covering, we know that
\begin{equation}
\label{step2eq1}
\I_{\O} (F) \textrm{ is close to } \sum_{j=1}^{\N}\alpha_{d-1}
r_j^{d-1} \nu(v_j) \,.
\end{equation}
On the other hand, thanks to property (ii), we obtain that
$$\E_n \cap B_j \,\approx\, F\cap B_j\, \approx\, B_j^-$$
for the distance $d$. It means
that in volume, $E_n$ is very similar to $B_j^-$ inside $B_j$, however
there might exist some thin but long strands in $B_j$ that belongs to $E_n
\cap (B_j^-)^c $ or to $E_n^c \cap B_j^-$. We want to compare $V(\E_n
\cap B_j)$ with the maximal flow $\tau_n (D_j, \gamma)$ in a cylinder
of basis $D_j = \disc (x_j, r_j', v_j)$ where $r_j'$ is a little bit
smaller than $r_j$, and $\gamma$ is a very small height, so
that the cylinder is included in $B_j$ and is almost flat. To make
this comparison, we have to cut the above-mentioned strands by adding edges to $\E_n$. We do it very carefully, in order to
control the number of edges we add, together with their capacity, and we
obtain that
\begin{equation}
\label{step2eq2}
V(\E_n \cap B_j) \,\leq\, \tau_n (D_j, \gamma) + \textrm{error}\,,
\end{equation}
where $error$ is a corrective term that is very small. Combining
(\ref{step2eq1}) and (\ref{step2eq2}), since $\I_{\O}(F)\geq \phi_{\O}$, we conclude that if $\phi_n
\leq (\phi_{\O} - \eps) n^{d-1} $ and $E_n \approx F$, then there exists
$j \in \{ 1,...,\N\}$ such that
$$\tau_n (D_j, \gamma) \,\leq\, (\nu(v_j) -
\eps/2) \alpha_{d-1} r_j'^{d-1} n^{d-1} \,.$$\\

{\bf Third step:} It remains to study the probability
$$ \PP [\tau_n (D_j, \gamma) \leq (\nu(v_j) -
\eps/2) \alpha_{d-1} r_j'^{d_1} n^{d-1}]\,. $$
In fact it has already been done by Rossignol and Th\'eret in
\cite{RossignolTheret08b}. It is easy to compare $\tau_n (D_j,
\gamma)$ with a sum of maximal flows through cylinders whose
bases are hyperrectangles. Then, we can use directly Theorem 3.9 in
\cite{RossignolTheret08b} that states that the lower large deviations
of these maximal flows below their limits are of surface order.


\subsubsection{Upper large deviations}
\label{secuppersketch}

To prove Theorem \ref{chapitre7devsup}, we have to study the
probability
\begin{equation}
\label{proba4}
 \PP \left[ \phi_n \geq (\widetilde{\phi_{\O}} + \eps) n^{d-1}\right]
\end{equation}
for a positive $\eps$. First of all, we can check that
$\widetilde{\phi_{\O}}$ is finite. In fact, we have to construct a
polyhedral set $P$ that satisfies all the conditions in
the definition of $\widetilde{\phi_{\O}}$. This is done with the help of
techniques very
similar to some of those we will use in section \ref{chapitre7secap}
to complete our polyhedral approximation, so we will not explain these
techniques here.
The proof of theorem \ref{chapitre7devsup} is divided in three steps.\\

{\bf First step:} We consider a polyhedral set $P$ as in the
definition of $\widetilde{\phi_{\O}}$ such that $\I_{\O} (P)$ is very
close to this constant. We want to construct sets of edges near
$\partial P \cap \O$ that cut $\Gamma^1_n$ from $\Gamma^2_n$ in
$\O_n$. Because we took a discrete approximation of $\O$ from the
outside, we need to enlarge a little $\O$, because some flow might go
from $\G^1_n$ to $\G^2_n$ using paths that lies partly in $\O_n
\smallsetminus \O$. Thus we construct a set $\O'$ which contains a small
neighbourhood of $\O$ (hence also $\O_n$ for all $n$ large
enough), which is transverse to $\partial P$, and which is small
enough to ensure that $\I_{\O'} (P)$ is still very close to $\phi_{\O}$. To
construct this set, we cover $\partial \O$ with small cubes,
by compactness we extract a finite subcover of $\partial \O$, and finally we
add the cubes of the subcover to $\O$ to obtain $\O'$. We construct these
cubes so that their boundaries are transverse to $\partial P$, and
their diameters are uniformly smaller than a small constant, so that
$\O'$ is included in a neighbourhood of $\O$ as small as we
need. Since $\partial P$ is transverse to $\G$, if we take this
constant small enough, we can control $\H^{d-1} (\partial P \cap (\O'
\smallsetminus \O))$, and thus the difference between $\I_{\O'} (P)$
and $\I_{\O} (P)$.

Then we construct a family of $C n$ (where $C>0$) disjoint sets of edges that
cut $\Gamma^1_n$ from $\Gamma^2_n$ in $\O_n$, and that lie near
$\partial P$. We consider the neighbourhood $P'$ of $P$ inside $\O'$
at distance smaller than a tiny constant $h$, and we partition $P' \smallsetminus P$ into
slabs $\M'(k)$ of width of order $1/n$, so we have $ Cn$ such slabs which look
like translates of $\partial P \cap \O'$ that are slightly deformed
and thickened. We prove that each
path from $\G^1_n$ to $\G^2_n$ in $\O_n$ must contain at least one
edge that lies in the set $\M'(k)$ for each $k$, i.e., each set
$\M'(k)$ contains a cutset. Thus we have found a family of $Cn$
disjoint cutsets.\\

{\bf Second step:} We almost cover $\partial P \cap \O'$ by a finite family
of disjoint cylinders $\B_j$, $j\in J$, whose bases are hyperrectangles
of sidelength $l$, that are orthogonal to $\partial P$,
of height bigger than $h$, and such that the part of $\partial P$
which is missing in this covering is very small. Thus, we obtain that
\begin{equation}
\label{step2beq1}
\I_{\O'} (P) \quad \textrm{is close to}\quad \sum_{j\in J} \nu(v_j) l^{d-1} \,,
\end{equation}
where $v_j$ gives the direction towards which the cylinder $B_j$ is
tilted (it is the unit vector which is orthogonal to the face of
$\partial P$ that cuts $B_j$).

We want to compare $\phi_n$ with the sum of the maximal flows $\phi
(B_j, v_j)$. For each $j$, let $E_j$ be a set of edges that cuts the
top from the bottom of $B_j$. The set $\cup_{j\in J} E_j$ does not cut
$\G^1_n$ from $\G^2_n$ in $\O_n$ in general, to create such a cutset
we must add two sets of edges:
\begin{itemize}
\item[(i)] a set of edges that covers the part of $\partial P \cap
  \O'$ that is missing in the covering by the cylinders $B_j$,
\item[(ii)] a set of edges that glues together all the previous sets of
  edges (the sets $E_j$ and the set described in (i)).
\end{itemize}
In fact, we have already constructed $Cn$ possible sets of edges as in (i): the edges that lie in
$\M'(k)\smallsetminus (\cup_{j\in J} B_j)$ for $k=1,...,Cn$. We
denote these sets by $M(k)$. We can also find $C' n$ ($C'>0$)
disjoint sets of edges that can be the glue described
in (ii), we denote these sets by $W(l)$ for $l=1,...,C'n$. We do not
provide a precise description of these sets. In fact, we can
choose different sets because we provide the glue more or
less in the interior of the cylinders $B_j$. Thus we obtain that
$$ \forall k\in \{1,...,Cn\}\,\, \forall l\in \{1,...,C'n\} \qquad
 \bigcup_{j\in J} E_j  \cup M(k) \cup W(l) \quad
\textrm{cuts }\G^1_n \textrm{ from } \G^2_n \textrm{ in }\O_n\,. $$
We obtain that
\begin{equation}
\label{step2beq2}
\phi_n \,\leq\, \sum_{j\in J} \phi(B_j, v_j) + \min_{k=1,...,Cn}
V(M(k)) + \min_{l=1,...,C'n} V(W(l))\,.
\end{equation}
Combining (\ref{step2beq1}) and (\ref{step2beq2}), we see that if
$\phi_n \geq (\widetilde{\phi_{\O}} +\eps ) n^{d-1}$, one of the
following events must happen:
\begin{itemize}
\item[(a)] $\exists j\in J \quad  \phi(B_j, v_j) \,\geq\, (\nu(v_j) +
  \eps/2) l^{d-1} n^{d-1} $,
\item[(b)] $\forall k\in\{1,...,Cn \} \quad V(M(k)) \,\geq\, \eta
  n^{d-1} $,
\item[(c)] $\forall l\in\{1,...,C'n \} \quad V(W(l)) \,\geq\, \eta
  n^{d-1} $,
\end{itemize}
where $\eta$ is a very small constant (depending on $\eps$ and $\phi_{\O}$).\\

{\bf Third step:} it consists in taking care of the
probability that the events (a), (b) or (c) happen. The probability of
(a) has already been studied in \cite{Theret:uppertau}: the upper
large deviations of the variable $\phi$ in a cylinder above $\nu$ are of volume order. The events (b) and (c) are of the same type,
and their probability is of the form
\begin{equation}
\label{step3beq1}
\PP \left[ \sum_{m=1}^{\alpha n^{d-1}} t_m \geq \eta n^{d-1}  \right]^{Dn}\,,
\end{equation}
where $(t_m)_{m\in \NN}$ is a family of i.i.d. variables of
  distribution function $\Lambda$, $D$ is a constant, $\eta$ is a very
  small constant and $\alpha n^{d-1}$ is the cardinality of the family
  of variables we consider. If $\alpha < \eta \EE[t_1]^{-1}$, and if the
  law $\Lambda$ admits one exponential moment, the
  Cram\'er Theorem in $\RR$ states that the probability
  (\ref{step3beq1}) decays exponentially fast with $n^d$. Note the
  role of the optimization over $Dn$ different probabilities to
  obtain the correct speed of decay. To complete the proof, it is enough to
  control the cardinality of the sets $M(k)$ and $W(l)$ for each $k$,
  $l$. This can been done, using the geometrical properties of
  $\partial P$ (it is polyhedral and transverse to $\partial \O'$).


\section[Polyhedral approximation]{Polyhedral approximation :
  $\phi_{\O} = \widetilde{\phi_{\O}}$}
\label{chapitre7secap}

We consider an open bounded domain $\O$ in $\RR^d$.
We denote its topological boundary by 
$\G=\partial\O$.
Let also
$\G^1$,
$\G^2$ 
be two
disjoint subsets of $\G$.
\medskip

\noindent
{\bf Hypothesis on $\O$:} 
We suppose that $\O$ is a Lipschitz domain, i.e., its boundary
$\G$ can be locally represented as the graph of a Lipschitz
function defined on some open ball of $\RR^{d-1}$.
Moreover there exists 
a finite number of oriented
hypersurfaces $S_1,\dots, S_p$ of class $C^1$
which are transverse to each other and
such that $\G$ is included in their union $S_1\cup\cdots\cup S_p$.
\medskip

\noindent
This
hypothesis is automatically satisfied when $\O$ is a 
bounded open set  with a $C^1$ boundary or when $\O$ is a polyhedral
domain. 
The Lipschitz condition can be expressed as follows:
each point $x$ of $\G=\partial \O$ has a neighbourhood
$U$ such that $U\cap\O$ is represented by the inequality
$x_n<f(x_1,\cdots,x_{n-1})$
in some cartesian coordinate
system where $f$ is a function satisfying a Lipschitz
condition.
Such domains are usually called Lipschitz domains in the
literature.
The boundary $\G$ of a Lipschitz domain is $d-1$ rectifiable
(in the terminology of Federer's book \cite{FED}),
so that its Minkowski content is equal to $\H^{d-1}(\G)$.
In addition, a Lipschitz domain $\O$ is admissible
(in the terminology of Ziemer's book \cite{ZI})
and in particular $\H^{d-1}(\G\smallsetminus\rb\O)=0$.
Moreover, each point of $\G$ is accessible from $\O$
through a rectifiable arc.  
\medskip

\noindent
{\bf Hypothesis on $\G^1,\G^2$:} 
The sets 
$\G^1$,
$\G^2$ 
are open subsets of $\G$. 
The relative boundaries 
$\partial_\G\,\,\G^1$,
$\partial_\G\,\,\G^2$ 
of $\G^1$, $\G^2$ 
in $\G$
have null $\H^{d-1}$ measure. The distance between $\G^1$ and $\G^2$
is positive.
\medskip

\noindent
We recall that the relative topology of $\G$ is
the topology induced on $\G$ by the topology
of $\RR^d$. Hence each of the sets  $\G^1,\G^2$ is the intersection
of $\G$ with an open set of $\RR^d$.
For $F$ a subset of $\O$ having finite perimeter in $\O$, 
the capacity of $F$ is
$$
\I_{\O}(F)\,=\,
\int_{\O\cap\rb F}
\kern-20pt
\nu(v_{F}(y))\,d\H^{d-1}(y)+
\int_{\G^2\cap\rb F}
\kern-20pt
\nu(v_F(y))\,d\H^{d-1}(y)+
\int_{\G^1\cap\rb(\O\smallsetminus F)}
\kern-20pt
\nu(v_{\O\smallsetminus F}(y))\,d\H^{d-1}(y)\,.
$$
For all $A\subset \RR^d$, $\overline{A}$ is the closure of $A$,
$\overset{\circ}{A}$ its interior and $A^c = \RR^d \smallsetminus A$.
We will prove the following theorem:
\begin{thm}\label{chapitre7rftf}
Let~$F$ be a subset of $\O$ having finite perimeter.
For any $\eps>0$, there exists a polyhedral set~$P$ whose boundary $\partial P$
is transverse to $\G$ and such that
$$ \overline{\G^1} \subset  \Po\ \,,\quad \overline{\G^2} \subset
\Pco \,,\quad \vol(F\Delta (P\cap \O)) \,<\,\eps\,,$$
$$\int_{\p^* P \cap \O} \nu(v_P(x)) d\H^{d-1}(x) \,=\, \I_{\O}(P)\,\leq\,\I_{\O}(F)+\eps \,.  $$
\end{thm}
\noindent

First we notice that theorem \ref{chapitre7rftf} implies that
$\phi_{\O} = \widetilde{\phi_{\O}}$, and thus the convergence of
$\phi_n$ (see section \ref{secres}). It is obvious since $\phi_{\O}
\leq \widetilde{\phi_{\O}}$ (see section \ref{seccontmincut}), and
theorem \ref{chapitre7rftf} implies that $\phi_{\O} \geq
\widetilde{\phi_{\O}}$.

The main difficulty of the proof of theorem \ref{chapitre7rftf} is to handle properly the
approximation close to $\G$ in order to push back inside $\O$
all the interfaces.
The essential tools of the proof are 
the Besicovitch differentiation theorem,
the Vitali covering theorem 
and an approximation technique due to De Giorgi.
Let us summarise the global strategy.
\medskip

\begin{figure}[!ht]
\centering
\begin{picture}(0,0)%
\includegraphics{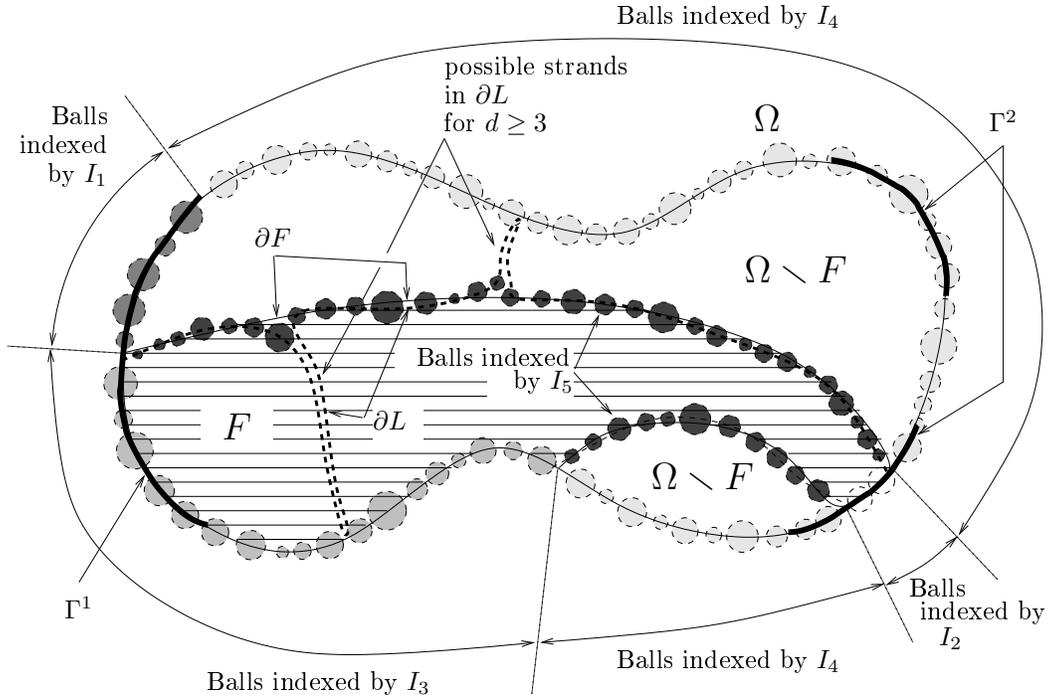}%
\end{picture}%
\setlength{\unitlength}{2368sp}%
\begingroup\makeatletter\ifx\SetFigFont\undefined%
\gdef\SetFigFont#1#2#3#4#5{%
  \reset@font\fontsize{#1}{#2pt}%
  \fontfamily{#3}\fontseries{#4}\fontshape{#5}%
  \selectfont}%
\fi\endgroup%
\begin{picture}(10858,7302)(814,-7582)
\put(10261,-6511){\makebox(0,0)[lb]{\smash{{\SetFigFont{10}{12.0}{\rmdefault}{\mddefault}{\updefault}{\color[rgb]{0,0,0}Balls}%
}}}}
\put(10576,-7021){\makebox(0,0)[lb]{\smash{{\SetFigFont{10}{12.0}{\rmdefault}{\mddefault}{\updefault}{\color[rgb]{0,0,0}$I_2$}%
}}}}
\put(10381,-6751){\makebox(0,0)[lb]{\smash{{\SetFigFont{10}{12.0}{\rmdefault}{\mddefault}{\updefault}{\color[rgb]{0,0,0}indexed by}%
}}}}
\put(3226,-4861){\makebox(0,0)[b]{\smash{{\SetFigFont{14}{16.8}{\rmdefault}{\mddefault}{\updefault}{\color[rgb]{0,0,0}$F$}%
}}}}
\put(3601,-2836){\makebox(0,0)[b]{\smash{{\SetFigFont{10}{12.0}{\rmdefault}{\mddefault}{\updefault}{\color[rgb]{0,0,0}$\p F$}%
}}}}
\put(5251,-7486){\makebox(0,0)[rb]{\smash{{\SetFigFont{10}{12.0}{\rmdefault}{\mddefault}{\updefault}{\color[rgb]{0,0,0}Balls indexed by $I_3$}%
}}}}
\put(1876,-1861){\makebox(0,0)[rb]{\smash{{\SetFigFont{10}{12.0}{\rmdefault}{\mddefault}{\updefault}{\color[rgb]{0,0,0}indexed}%
}}}}
\put(1876,-2161){\makebox(0,0)[rb]{\smash{{\SetFigFont{10}{12.0}{\rmdefault}{\mddefault}{\updefault}{\color[rgb]{0,0,0}by $I_1$}%
}}}}
\put(9076,-3211){\makebox(0,0)[b]{\smash{{\SetFigFont{14}{16.8}{\rmdefault}{\mddefault}{\updefault}{\color[rgb]{0,0,0}$\Omega \smallsetminus F$}%
}}}}
\put(7201,-511){\makebox(0,0)[lb]{\smash{{\SetFigFont{10}{12.0}{\rmdefault}{\mddefault}{\updefault}{\color[rgb]{0,0,0}Balls indexed by $I_4$}%
}}}}
\put(7201,-7261){\makebox(0,0)[lb]{\smash{{\SetFigFont{10}{12.0}{\rmdefault}{\mddefault}{\updefault}{\color[rgb]{0,0,0}Balls indexed by $I_4$}%
}}}}
\put(8101,-5386){\makebox(0,0)[b]{\smash{{\SetFigFont{14}{16.8}{\rmdefault}{\mddefault}{\updefault}{\color[rgb]{0,0,0}$\Omega \smallsetminus F$}%
}}}}
\put(8776,-1636){\makebox(0,0)[b]{\smash{{\SetFigFont{14}{16.8}{\rmdefault}{\mddefault}{\updefault}{\color[rgb]{0,0,0}$\Omega$}%
}}}}
\put(1576,-6736){\makebox(0,0)[b]{\smash{{\SetFigFont{10}{12.0}{\rmdefault}{\mddefault}{\updefault}{\color[rgb]{0,0,0}$\Gamma^1$}%
}}}}
\put(11251,-1636){\makebox(0,0)[b]{\smash{{\SetFigFont{10}{12.0}{\rmdefault}{\mddefault}{\updefault}{\color[rgb]{0,0,0}$\Gamma^2$}%
}}}}
\put(5401,-1636){\makebox(0,0)[lb]{\smash{{\SetFigFont{10}{12.0}{\rmdefault}{\mddefault}{\updefault}{\color[rgb]{0,0,0}for $d\geq 3$}%
}}}}
\put(5401,-1336){\makebox(0,0)[lb]{\smash{{\SetFigFont{10}{12.0}{\rmdefault}{\mddefault}{\updefault}{\color[rgb]{0,0,0}in $\p L$}%
}}}}
\put(5401,-1036){\makebox(0,0)[lb]{\smash{{\SetFigFont{10}{12.0}{\rmdefault}{\mddefault}{\updefault}{\color[rgb]{0,0,0}possible strands}%
}}}}
\put(4651,-4756){\makebox(0,0)[lb]{\smash{{\SetFigFont{10}{12.0}{\rmdefault}{\mddefault}{\updefault}{\color[rgb]{0,0,0}$\p L$}%
}}}}
\put(6721,-4096){\makebox(0,0)[rb]{\smash{{\SetFigFont{10}{12.0}{\rmdefault}{\mddefault}{\updefault}{\color[rgb]{0,0,0}Balls indexed}%
}}}}
\put(1876,-1561){\makebox(0,0)[rb]{\smash{{\SetFigFont{10}{12.0}{\rmdefault}{\mddefault}{\updefault}{\color[rgb]{0,0,0}Balls}%
}}}}
\put(6751,-4336){\makebox(0,0)[rb]{\smash{{\SetFigFont{10}{12.0}{\rmdefault}{\mddefault}{\updefault}{\color[rgb]{0,0,0}by $I_5$}%
}}}}
\end{picture}%
\caption{The balls indexed by $I_i$ for $i=1,...,5$.}
\label{chapitre7ap}
\end{figure}

{\noindent {\bf Sketch of the proof:\ }}
We fix $\ga>0$.
We cover $\rb\O$ up to a set of $\H^{d-1}$ measure less than~$\ga$
by a finite collection of disjoint balls
$B(x_i,r_i)$, $i\in I_1\cup I_2\cup I_3\cup I_4$, 
centered on $\G$, whose radii are sufficiently small to ensure that 
the surface and volume estimates
within the balls are controlled by the factor $\ga$.
The indices of $I_1$ correspond to
balls centered on $\G^1\cap \rb (\O\smallsetminus F)$,
the indices of $I_2$ to
balls centered on $\G^2\cap \rb F$,
the indices of $I_3$ to
balls centered on $(\G \smallsetminus \G^2 )\cap\rb F $,
the indices of $I_4$ to
balls centered on $(\G \smallsetminus \G^1 )\cap\rb (\O\smallsetminus F)$
(see figure \ref{chapitre7ap}).
The remaining part of $\G$ is covered 
by a finite collection of balls
$B(y_j,s_j)$, $j\in J_0\cup J_1\cup J_2$.
The indices of $J_1$ correspond to
balls covering the remaining part of $\overline{\G}_1$,
the indices of $J_2$ correspond to
balls covering the remaining part of $\overline{\G}_2$.
We choose $\eps>0$ sufficiently small, depending on $\ga$ 
and on the previous families of balls and
we approximate
the set $F$ by a smooth set $L$ inside $\O$, whose capacity and
volume are at distance less than $\eps$ from those of $F$.
We build then two further family of balls:
\begin{itemize}
\item[-] $B(x_i,r_i)$, $i\in I_5$, cover $\O\cap\partial L$,
up to a set of $\H^{d-1}$ measure $\eps$.
\item[-] $B(y_j,s_j)$, $j\in J_3$, cover the remaining 
set
$\O\cap\partial L\smallsetminus
\bigcup_{i\in I_5}B(x_i,r_i)$.
\end{itemize}
Inside each ball 
$B(x_i,r_i)$, 
$i\in I_1\cup I_2 \cup I_3\cup I_4\cup I_5$,
up to a small fraction, the interfaces
are located on hypersurfaces and the radii of the balls
are so small that these hypersurfaces are almost flat.
Hence we can enclose the interfaces into small flat
polyhedral cylinders 
$D_i$, 
$i\in I_1\cup I_2 \cup I_3\cup I_4\cup I_5$,
and by aggregating adequately the cylinders to the
set $F$ or to its complement $\O\smallsetminus F$, 
we move these interfaces on the
boundaries of these cylinders.
The remaining interfaces are enclosed in the balls
$B(y_j,s_j)$, $j\in J_0\cup J_1\cup J_2\cup J_3$ and we 
approximate these balls from the outside
by polyhedra.

We have to define delicately the whole process, in order not to lose too
much capacity, and to control the possible interaction between
interfaces close to $\G$ and interfaces in $\O$.
The presence of boundary conditions creates a substantial additional
difficulty
compared to the polyhedral approximation
performed in \cite{Cerf:StFlour}.
Indeed, the most difficult interfaces to handle are those
corresponding to $D_i$, $i\in I_3\cup I_4$. We first choose the
balls $B(x_i,r_i)$, $i\in I_1\cup I_2\cup I_3\cup I_4$, 
corresponding to $\ga$.
We cover the remaining portion of~$\G$ with the balls
$B(y_j,s_j)$, $j\in J_0\cup J_1\cup J_2$.
At this point we can already in principle
define the cylinders
$D_i$, $i\in I_1\cup I_2$. Then we choose $\eps$ small enough,
depending on $\ga$ and the balls
$B(x_i,r_i)$, $i\in I_1\cup I_2\cup I_3\cup I_4$, 
to ensure that the perturbation of volume $\eps$ caused 
when smoothing the set $F$ inside $\O$
will not alter significantly
the situation inside the balls
$B(x_i,r_i)$, $i\in I_3\cup I_4$.
Then we move inside $\O$ and we build the cylinders
$D_i$, $i\in I_5$. Then we come back to the boundary and
we build the cylinders
$D_i$, $i\in I_3\cup I_4$. We cover the remaining interfaces 
in $\O$ by the balls
$B(y_j,s_j)$, $j\in J_3$.
Finally we aggregate successively each flat polyhedral cylinder
$D_i$ to the set $L$ or to its complement.
\par
\medskip

{\noindent {\bf Preparation of the proof.}}
Let us consider a subset $F$ of $\O$ having finite perimeter.
Let $\gamma$ belong to $]0,1/16[$.
We start by handling the boundary $\G$, for which we make
locally flat approximations controlled by the factor $\ga$.
By hypothesis, there exists
a finite number of oriented
hypersurfaces $S_1,\dots, S_p$ of class $C^1$
such that $\G$ is included in their union $S_1\cup\cdots\cup S_p$.
In particular, we have
$$\G\smallsetminus\rb\O\,\subset\,S\,=\,
\bigcup_{1\leq k<l\leq p}S_k\cap S_l\,.$$
Since the hypersurfaces $S_1,\dots,S_r$ are transverse to each other,
this implies that
$\H^{d-1}(S)=0$.

{\noindent {\bf $\bullet$ Continuity of the normal vectors.}}
The hypersurfaces $S_1,\dots,S_p$ being $C^1$ and the set $\G$
compact, the maps $x\in\G\mapsto v_{S_k}(x)$, $1\leq k\leq p$
(where $v_{S_k}(x)$ is the unit normal vector to $S_k$ at $x$)
are uniformly continuous:
$$\forall \delta>0\quad\exists\eta>0\quad
\forall 
k\in\{\,1,\dots,p\,\}
\quad
\forall x,y\in S_k\cap\G\quad
|x-y|_2\leq \eta\,\Rightarrow\,
\big|v_{S_k}(x)-v_{S_k}(y)\big|_2<\delta
\,.$$
Let $\eta^*$ be associated to $\delta=1$ by this property.
We will use also a more refined property.

{\noindent {\bf $\bullet$ Localisation of the interfaces.}}
We first prove a geometric lemma:
\begin{lem}\label{chapitre7rfpres}
Let~$\Gamma$ be an hypersurface
(that is a $C^1$ submanifold of~$\RR^d$ of codimension~$1$)
and let~$K$ be a compact subset of~$\Gamma$.
There exists a positive~$M=M(\Gamma,K)$ such that:
$$\forall\eps>0\quad\exists\,r>0\quad
\forall x,y\in K\qquad
|x-y|_2\leq r\quad\Rightarrow\quad
d_2(y,\tan(\Gamma,x))\leq M\,\eps\,|x-y|_2\,.$$
($\tan(\Gamma,x)$ is the tangent hyperplane of~$\Gamma$ at~$x$).
\end{lem}
\proof
By a standard compactness argument, it is enough to prove
the following local property:
$$\displaylines{
\forall x\in\Gamma\quad
\exists\,M(x)>0\quad\,\forall\eps>0\quad
\exists\,r(x,\eps)>0\quad
\forall y,z\in \Gamma\cap B(x,r(x,\eps))\qquad
\hfill\cr
\hfill
d_2(y,\tan(\Gamma,z))\leq M(x)\,\eps\,|y-z|_2\,.}$$
Indeed, if this property holds, we cover $K$ 
by the open balls $\bro (x,r(x,\eps)/2)$, $x\in K$,
we extract a finite subcover
$\bro(x_i,r(x_i,\eps)/2)$, $1\leq i\leq k$,
and we set
$$M=\max\{\,M(x_i):1\leq i\leq k\,\}\,,\quad
r=\min\{\,r(x_i,\eps)/2:1\leq i\leq k\,\}\,.$$
Let now~$y,z$ belong to~$K$ with $|y-z|_2\leq r$.
Let~$i$ be such that $y$ belongs to
$B(x_i,r(x_i,\eps)/2)$.
Since $r \leq r(x_i,\eps)/2$,
then both~$y,z$ belong to the ball
$B(x_i,r(x_i,\eps))$ and it follows that
$$d_2(y,\tan(\Gamma,z))\,\leq\, M(x_i)\,\eps\,|y-z|_2\leq
M\,\eps\,|y-z|_2\,.$$
\par

We turn now to the proof of the above local property.
Since~$\Gamma$ is an hypersurface, for any~$x$ in~$\Gamma$ there
exists a neighbourhood~$V$ of~$x$ in~$\RR^d$, a diffeomorphism
$f:V\mapsto\RR^d$ of class~$C^1$ and a $(d-1)$~dimensional 
vector space~$Z$ of~$\RR^d$ such that
$Z\cap f(V)=f(\Gamma\cap V)$ (see for instance \cite{FED}, $3.1.19$).
Let~$A$ be a compact neighbourhood of~$x$ included in~$V$.
Since $f$ is a diffeomorphism, the maps
$y\in A\mapsto df(y)\in\text{End}(\RR^d)$,
$u\in f(A)\mapsto df^{-1}(u)\in\text{End}(\RR^d)$
are continuous. 
Therefore they are bounded:
$$\exists M>0\quad\forall y\in A\quad ||df(y)||\leq M\,,\quad
\forall u\in f(A)\quad ||df^{-1}(u)||\leq M\,$$
(here 
$||df(x)||=\sup\{\,|df(x)(y)|_2:|y|_2\leq 1\,\}$
is the standard operator norm in $\text{End}(\RR^d)$).
Since $f(A)$ is compact,
the differential map $df^{-1}$ is uniformly continuous on~$f(A)$:
$$\forall\eps>0\quad\exists\delta>0\quad
\forall u,v\in f(A)\qquad
|u-v|_2\leq\delta\quad\Rightarrow\quad
||df^{-1}(u)-df^{-1}(v)||\leq\eps\,.$$
Let~$\eps$ be positive and let~$\delta$ be associated to~$\eps$
as above.
Let~$\rho$ be positive and small enough so that
$\rho<\delta/2$ and $B(f(x),\rho)\subset f(A)$
(since $f$ is a $C^1$ diffeomorphism, $f(A)$ is a 
neighbourhood of $f(x)$).
Let~$r$ be such that $0<r<\rho/M$ and
$B(x,r)\subset A$. We claim that~$M$ associated to~$x$ and
$r$ associated to~$\eps,x$ answer the problem.
Let~$y,z$ belong to~$\Gamma\cap B(x,r)$.
Since $[y,z]\subset B(x,r)\subset A$, 
and $||df(\zeta)||\leq M$ on $A$, then
$$\displaylines{
|f(y)-f(x)|_2\leq M|y-x|_2\leq Mr<\rho\,,\quad
|f(z)-f(x)|_2<\rho\,,\quad\cr
|f(y)-f(z)|_2<\delta\,,\quad
|f(y)-f(z)|_2<M|y-z|_2\,.}$$
We apply next a classical lemma 
of differential calculus (see \cite{LA}, I, 4, Corollary~$2$)
to the map $f^{-1}$ and the
interval $[f(z),f(y)]$ (which is included in $B(f(x),\rho)\subset f(A)$)
and the point $f(z)$:
$$\displaylines{
|y-z-df^{-1}(f(z))(f(y)-f(z))|_2\,\leq\,\hfill\cr
\hfill
|f(y)-f(z)|_2
\sup\,\{\,||df^{-1}(\zeta)-df^{-1}(f(z))||:\zeta\in [f(z),f(y)]\,\}\,.}$$
The right--hand member is less than
$M|y-z|_2\,\eps$.
Since 
$z+df^{-1}(f(z))(f(y)-f(z))$
belongs to $\tan(\Gamma,z)$, we are done.\qed

We come back to our case.
Let $k\in\{\,1,\dots,p\,\}$. The set $S_k\cap\G$ is a compact subset
of the hypersurface $S_k$.
Applying lemma~\ref{chapitre7rfpres}, we get:
$$\displaylines{
\exists M_k\,\,
\forall\delta_0>0\,\, \exists\,\eta_k>0\,\,
\forall x,y\in 
S_k\cap\G
\quad
|x-y|_2\leq \eta_k\,\Rightarrow\,
d_2\big(y,\tan(S_k,x)\big)
\leq M_k\delta_0 |x-y|_2\,.}$$
Let $M_0\,=\,\max_{1\leq k\leq p}M_k$ and 
let $\delta_0$ in $]0,1/2[$
be such that $M_0\delta_0<\gamma$.
For each $k$ in $\{\,1,\dots,p\,\}$, let 
$\eta_k$ be associated to
$\delta_0$ as in the above property and let
$$\eta_0\,=\,\min\Big(\min_{1\leq k\leq p}\eta_k,
\,\eta^*,\,
\frac{1}{8d}\,\text{dist}
(\G^1,\G^2)
\Big)\,.$$
{\noindent {\bf $\bullet$ Covering of $\G$ by transverse cubes.}}
We build a family of cubes $Q(x,r)$, indexed by $x\in \G$
and $r\in ]0,r_\G[$ such that
$Q(x,r)$ is a cube centered at $x$ of side length $r$ which is transverse
to~$\G$.
For $x\in \RR^d$ and $k\in\{\,1,\dots,p\,\}$, let $p_k(x)$ be a
point of $S_k\cap\G$ such that
$$|x-p_k(x)|_2\,=\,
\inf\,\big\{\,|x-y|_2:y\in S_k\cap\G\,\big\}\,.$$
Such a point exists since 
$S_k\cap\G$ is compact.
We define then for
$k\in\{\,1,\dots,p\,\}$
$$\forall x\in \RR^d\qquad
v_k(x)\,=\,v_{S_k}(p_k(x))\,.$$
We define also
$$d_r\,=\,
\inf_{v_1,\dots,v_p\in S^{d-1}}\,
\max_{b\in{\mathcal B}_d}\,\,\,
\min_{\begin{array}{c} 1\leq k\leq r \\  e\in b \end{array}}\,
\big(|e-v_i|_2,|-e-v_i|_2\big)$$
where  
${\mathcal B}_d$
is the collection of the orthonormal basis of ${\mathbb R}^d$
and $S^{d-1}$ is the 
unit sphere of $\RR^d$.
Let $\eta$ be associated to $d_r/4$ as in the 
above continuity property. We set 
$$r_\G=\frac{\eta}{2d}\,.$$
Let $x\in \G$. By the definition of $d_r$, there exists 
an orthonormal basis $b_x$ of $\RR^d$ such that
$$\forall e\in b_x\quad
\forall k\in\{\,1,\dots,p\,\}
\quad
\min\big(|e-v_k(x)|_2,
|-e-v_k(x)|_2\big)\,>\,\frac{d_r}{2}\,.$$
Let $Q(x,r)$ be the cube centered at $x$ of sidelength~$r$
whose sides are parallel to the vectors of $b_x$.
We claim that $Q(x,r)$ is transverse to~$\G$ for $r<r_\G$.
Indeed, let $y\in Q(x,r)\cap\G$. Suppose that $y\in S_k$ for some
$k\in\{\,1,\dots,p\,\}$, so that
$v_k(y)=v_{S_k}(y)$ and 
$|x-p_k(x)|_2<dr_\G$. In particular, we have
$|y-p_k(x)|_2<2dr_\G<\eta$ and
$|v_{S_k}(y)-v_k(x)|_2<d_r/4$. 
For $e\in b_x$,
$$
\frac{d_r}{2}\,\leq\,
|e-v_k(x)|_2
\,\leq\,
|e-v_{S_k}(y)|_2+
|v_{S_k}(y)-v_k(x)|_2
$$
whence
$$|e-v_{S_k}(y)|_2
\,\geq\,
\frac{d_r}{2}-
\frac{d_r}{4}\,=\,
\frac{d_r}{4}
\,.$$
This is also true for $-e$, therefore
the faces of the cube $Q(x,r)$ are transverse to $S_k$.

{\noindent {\bf $\bullet$ Vitali covering Theorem for $\H^{d-1}$.}}
A collection of sets $\U$ is called a Vitali
class for a Borel set $E$ of $\RR^d$ if for each $x\in
E$ and $\delta >0$, there exists a set $U\in \U$ containing $x$ such that
$0<\diam U < \delta$, where $\diam U$ is the diameter of the set $U$. We
now recall the Vitali covering Theorem for $\H^{d-1}$ (see for
instance \cite{FAL}, Theorem 1.10), since it will be useful during
the proof:
\begin{thm}
\label{chapitre7vitali}
Let $E$ be a $\H^{d-1}$ measurable subset of $\RR^d$ and $\U$ be a Vitali
class of closed sets for $E$. Then we may select a (countable) disjoint
sequence $(U_i)_{i\in I}$ from $\U$ such that 
$$ \textrm{either } \sum_{i\in I} (\diam U_i)^{d-1} \,=\, +\infty \textrm{ or
} \H^{d-1} (E\smallsetminus \cup_{i\in I} U_i) \,=\, 0 \,.$$
If $\H^{d-1} (E) <\infty$, then given $\eps >0$, we may also require that
$$ \H^{d-1} (E) \,\leq\, \frac{\alpha_{d-1}}{2^{d-1}} \sum_{i\in I} (\diam
U_i)^{d-1} \,. $$
\end{thm}

\begin{figure}
\centering
\begin{picture}(0,0)%
\includegraphics{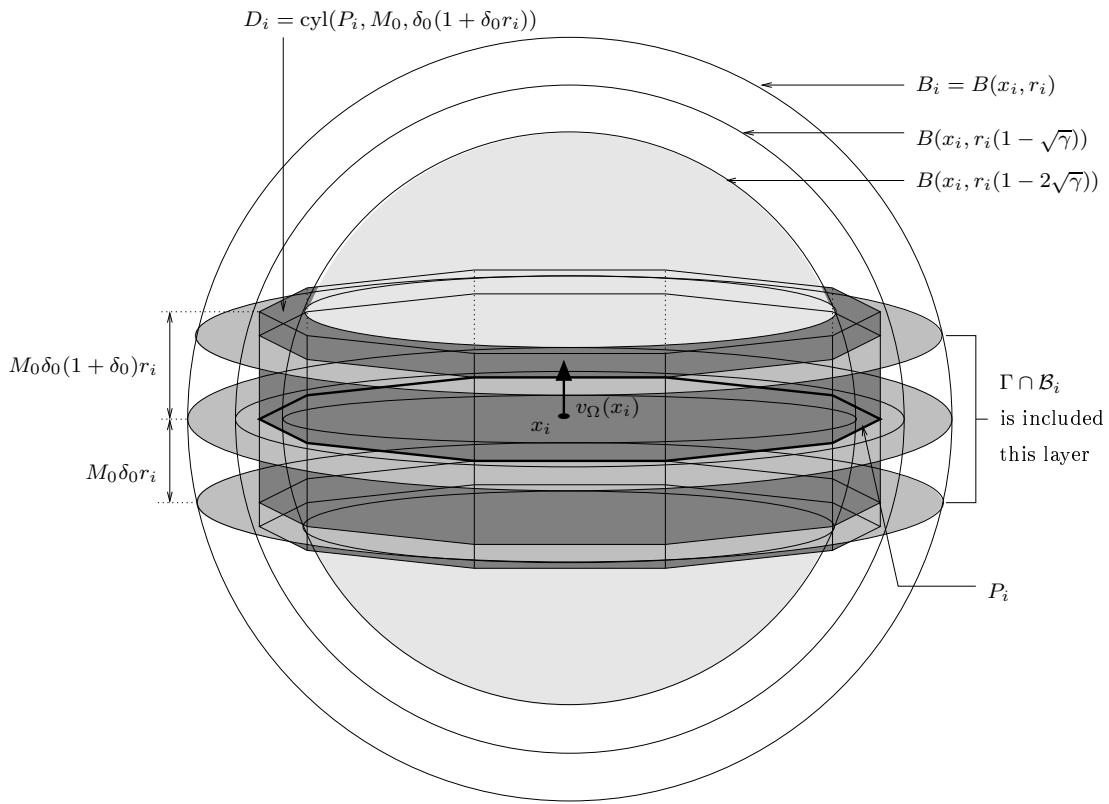}%
\end{picture}%
\setlength{\unitlength}{1973sp}%
\begingroup\makeatletter\ifx\SetFigFont\undefined%
\gdef\SetFigFont#1#2#3#4#5{%
  \reset@font\fontsize{#1}{#2pt}%
  \fontfamily{#3}\fontseries{#4}\fontshape{#5}%
  \selectfont}%
\fi\endgroup%
\begin{picture}(11820,9988)(811,-9968)
\put(11401,-4786){\makebox(0,0)[lb]{\smash{{\SetFigFont{8}{9.6}{\rmdefault}{\mddefault}{\updefault}{\color[rgb]{0,0,0}$\G \cap \B_i$}%
}}}}
\put(11401,-5236){\makebox(0,0)[lb]{\smash{{\SetFigFont{8}{9.6}{\rmdefault}{\mddefault}{\updefault}{\color[rgb]{0,0,0}is included}%
}}}}
\put(11401,-5686){\makebox(0,0)[lb]{\smash{{\SetFigFont{8}{9.6}{\rmdefault}{\mddefault}{\updefault}{\color[rgb]{0,0,0}this layer}%
}}}}
\put(11251,-7411){\makebox(0,0)[lb]{\smash{{\SetFigFont{8}{9.6}{\rmdefault}{\mddefault}{\updefault}{\color[rgb]{0,0,0}$P_i$}%
}}}}
\put(10351,-2236){\makebox(0,0)[lb]{\smash{{\SetFigFont{8}{9.6}{\rmdefault}{\mddefault}{\updefault}{\color[rgb]{0,0,0}$B(x_i, r_i (1-2\sqrt{\gamma}))$}%
}}}}
\put(10351,-1711){\makebox(0,0)[lb]{\smash{{\SetFigFont{8}{9.6}{\rmdefault}{\mddefault}{\updefault}{\color[rgb]{0,0,0}$B(x_i, r_i (1-\sqrt{\gamma}))$}%
}}}}
\put(10351,-1036){\makebox(0,0)[lb]{\smash{{\SetFigFont{8}{9.6}{\rmdefault}{\mddefault}{\updefault}{\color[rgb]{0,0,0}$B_i = B(x_i,r_i)$}%
}}}}
\put(826,-4561){\makebox(0,0)[rb]{\smash{{\SetFigFont{8}{9.6}{\rmdefault}{\mddefault}{\updefault}{\color[rgb]{0,0,0}$M_0 \delta_0 (1+\delta_0 )r_i$}%
}}}}
\put(826,-5911){\makebox(0,0)[rb]{\smash{{\SetFigFont{8}{9.6}{\rmdefault}{\mddefault}{\updefault}{\color[rgb]{0,0,0}$M_0 \delta_0 r_i$}%
}}}}
\put(3751,-211){\makebox(0,0)[b]{\smash{{\SetFigFont{8}{9.6}{\rmdefault}{\mddefault}{\updefault}{\color[rgb]{0,0,0}$D_i = \cyl(P_i, M_0, \delta_0 (1+\delta_0 r_i))$}%
}}}}
\put(5776,-5311){\makebox(0,0)[rb]{\smash{{\SetFigFont{8}{9.6}{\rmdefault}{\mddefault}{\updefault}{\color[rgb]{0,0,0}$x_i$}%
}}}}
\put(6076,-5086){\makebox(0,0)[lb]{\smash{{\SetFigFont{8}{9.6}{\rmdefault}{\mddefault}{\updefault}{\color[rgb]{0,0,0}$v_{\O}(x_i)$}%
}}}}
\end{picture}%
\caption{The cylinder $D_i$ for $i\in I_1 \cup I_2$.}
\label{chapitre7polygone3}
\end{figure}

{\noindent 
{\bf Start of the main argument}.}
We first handle the interfaces along~$\G$.
Let ${\R}(\G)$ be the set of the points $x$
of $\G\smallsetminus S$ such that
$$\displaylines{
\lim_{r\to 0}\quad
(\al_{d}r^{d})^{-1}
\ld(B(x,r)\smallsetminus\O)\,=\,1/2\,,\cr
\lim_{r\to 0}\quad
(\al_{d-1}r^{d-1})^{-1}
\H^{d-1}(B(x,r)\cap\G)\,=\,1\,.\cr
}$$
Let ${\R}(\O\smallsetminus F)$ be the set of the points $x$
belonging to $\rb(\O\smallsetminus F)\cap {\R}(\G)$
such that
$$\displaylines{
\lim_{r\to 0}\quad
(\al_{d-1}r^{d-1})^{-1}
\H^{d-1}(B(x,r)\cap\rb(\O\smallsetminus F))\,=\,1\,,\cr
\lim_{r\to 0}\quad
(\al_{d}r^{d})^{-1}
\ld(B(x,r)\cap (\O\smallsetminus F))\,=\,1/2\,,\cr
\lim_{r\to 0}\quad
(\al_{d-1}r^{d-1})^{-1}
\int_{B(x,r)\cap\rb(\O\smallsetminus F)}
\kern-20pt
\nu(v_{\O\smallsetminus F}(y))\,d\H^{d-1}(y)
\,=\,\nu(v_\O(x))\,.}$$
Let ${\R}(F)$ be the set of the points $x$
belonging to $\rb F\cap{\R}(\G)$
such that
$$\displaylines{
\lim_{r\to 0}\quad
(\al_{d-1}r^{d-1})^{-1}
\H^{d-1}(B(x,r)\cap\rb F)\,=\,1\,,\cr
\lim_{r\to 0}\quad
(\al_{d}r^{d})^{-1}
\ld(B(x,r)\cap F)\,=\,1/2\,,\cr
\lim_{r\to 0}\quad
(\al_{d-1}r^{d-1})^{-1}
\int_{B(x,r)\cap\rb F}
\kern-20pt
\nu(v_F(y))\,d\H^{d-1}(y)
\,=\,\nu(v_\O(x))\,.}$$
Thanks to the hypothesis on~$\G$ and
the structure of the sets of finite perimeter
(see either Lemma~$1$, section~$5.8$ of \cite{EVGA},
Lemma~$5.9.5$ in \cite{ZI} or
Theorem~$3.61$ of \cite{AMBRO}), 
we have 
$$\H^{d-1}\big(\G\smallsetminus 
({\R}(F)\cup {\R}(\O\smallsetminus F))\big)\,=\,0\,.$$
For~$x$ in~${\R}(\G)$,
there exists a positive~$r_0(x,\gamma)$ such that, 
for any $r<r_0(x,\gamma)$, 
$$\displaylines{
|\ld(B(x,r)\smallsetminus\O )-\al_{d}r^{d}/2|
\,\leq\,\gamma\,\al_{d}r^{d}\,,
\cr
|\H^{d-1}(B(x,r)\cap\G)-\al_{d-1}r^{d-1}|
\,\leq\,\gamma\,\al_{d-1}r^{d-1}\,.
}$$
For~$x$ in~${\R}(\O\smallsetminus F)$,
there exists a positive~$r(x,\gamma)<r_0(x,\ga)$ such that, 
for any $r<r(x,\gamma)$, 
$$\displaylines{
|\H^{d-1}(B(x,r)\cap\rb(\O\smallsetminus F))
-\al_{d-1}r^{d-1}|
\,\leq\,\gamma\,\al_{d-1}r^{d-1}
\,,\cr
|\ld(B(x,r)\cap(\O\smallsetminus F) )-\al_{d}r^{d}/2|
\,\leq\,\gamma\,\al_{d}r^{d}\,,\cr
\Big|
(\al_{d-1}r^{d-1})^{-1}
\int_{B(x,r)\cap\rb(\O\smallsetminus F)}
\kern-20pt
\nu(v_{\O\smallsetminus F}(y))\,d\H^{d-1}(y)
-\nu(v_\O(x))
\Big|\,\leq\,
\gamma\,
\,.}$$
For~$x$ in~${\R}(F)$,
there exists a positive~$r(x,\gamma)<r_0(x,\ga)$ such that, 
for any $r<r(x,\gamma)$, 
$$\displaylines{
|\H^{d-1}(B(x,r)\cap\rb F)
-\al_{d-1}r^{d-1}|
\,\leq\,\gamma\,\al_{d-1}r^{d-1}
\,,\cr
|\ld(B(x,r)\cap F )-\al_{d}r^{d}/2|
\,\leq\,\gamma\,\al_{d}r^{d}\,,\cr
\Big|
(\al_{d-1}r^{d-1})^{-1}
\int_{B(x,r)\cap\rb F}
\kern-20pt
\nu(v_{F}(y))\,d\H^{d-1}(y)
-\nu(v_\O(x))
\Big|\,\leq\,
\gamma\,
\,.}$$
Let us define the sets
$$\displaylines{
\G^{1*}\,=\,\G^1\cap{\R}(\O\smallsetminus F)\,,\quad
\G^{2*}\,=\,\G^2\cap{\R}(F)\,,\cr
\G^{3*}\,=\,(\G\smallsetminus\bGa_2)\cap{\R}(F)\,,\quad
\G^{4*}\,=\,(\G\smallsetminus\bGa_1)\cap{\R}(\O\smallsetminus F)\,.}$$
The family of balls
$$\displaylines{
B(x,r)\,,\quad
x\in\G^{1*}\cup\G^{2*}\,,\quad
r<\min\big(r(x,\gamma),\gamma,\eta_0,\frac{1}{2}\text{dist}(x,S)\big)\,,\cr
B(x,r)\,,\quad
x\in\G^{3*}\,,\quad
r<\min\big(r(x,\gamma),\gamma,\eta_0,\frac{1}{2}\text{dist}(x,S),
\frac{1}{ 2}\text{dist}(x,\bGa_2)
\big)\,,\cr
B(x,r)\,,\quad
x\in\G^{4*}\,,\quad
r<\min\big(r(x,\gamma),\gamma,\eta_0,\frac{1}{2}\text{dist}(x,S),
\frac{1}{ 2}\text{dist}(x,\bGa_1)
\big)\,\phantom{,}
}$$
is a Vitali relation for 
$\G^{1*}\cup\G^{2*}\cup\G^{3*}\cup\G^{4*}$. Recall that $S$ is the set
of the points belonging to two or more of 
the hypersurfaces $S_1,\dots,S_p$ and since $S$ is disjoint from
$\G^{1*},\G^{2*},\G^{3*},\G^{4*}$,
then $\text{dist}(x,S)>0$ for
$x\in\G^{1*}\cup\G^{2*}\cup\G^{3*}\cup\G^{4*}$.
By the standard Vitali covering Theorem (see theorem \ref{chapitre7vitali}),
we may select a finite or countable collection
of disjoint balls $B(x_i,r_i)$, $i\in I$, such that:
for $i\in I$, $x_i\in
\G^{1*}\cup\G^{2*}\cup\G^{3*}\cup\G^{4*}$,
$r_i<\min(r(x_i,\gamma),\gamma,\eta_0,
\frac{1}{2}\text{dist}(x_i,S))$
and
$$\text{either}\qquad\H^{d-1}\Big(
\G
\smallsetminus\bigcup_{i\in I}B(x_i,r_i)
\Big)\,=\,0\qquad\text{or}\qquad
\sum_{i\in I}r_i^{d-1}\,=\,\infty
\,.$$
Because for each~$i$ in~$I$, $r_i$ is smaller than 
$r(x_i,\gamma)$, 
$$\al_{d-1}(1-\gamma)\sum_{i\in I}r_i^{d-1}\,
\leq\,\H^{d-1}(\G)\,<\,\infty\,$$
and therefore the first case occurs, so that we may select 
four finite
subsets~$I_1,I_2, I_3, I_4$ of~$I$ such that
 $$\displaylines{
\forall k\in \{\,1,\dots,4\,\}\quad
\forall i\in I_k\quad x_i\in \G^{k*}\,,
\cr
\H^{d-1}\Big(
\G
\smallsetminus\bigcup_{1\leq k\leq 4}\bigcup_{i\in I_k}B(x_i,r_i)
\Big)\,<\,\gamma\,.}$$
Let $i$ belong to $I_1\cup I_2\cup I_3\cup I_4$. We have
$$\displaylines{
\H^{d-1}(\G\cap B(x_i,r_i)\smallsetminus
B(x_i,r_i(1-2\sqrt\ga)))=
\H^{d-1}(\G\cap B(x_i,r_i))-
\H^{d-1}(\G\cap B(x_i,r_i(1-2\sqrt\ga)))\cr
\leq (1+\ga)\al_{d-1}r_i^{d-1}-
(1-\ga)\al_{d-1}r_i^{d-1}(1-2\sqrt\ga)^{d-1}\cr
=\al_{d-1}r_i^{d-1}(1+\ga-(1-\ga)(1-2\sqrt\ga)^{d-1})\cr
\leq \al_{d-1}r_i^{d-1}2d\sqrt\ga\,.}$$
Hence
\begin{align*}
\sum_{i\in I_1\cup I_2\cup I_3\cup I_4}
\H^{d-1}(\G\cap B(x_i,r_i) & \smallsetminus
B(x_i,r_i(1-2\sqrt\ga)))\\ &\,\leq\,
2d\sqrt\ga\kern-3pt
\sum_{i\in I_1\cup I_2\cup I_3\cup I_4}
\kern-7pt
\al_{d-1}r_i^{d-1}
\,\leq\,
4d\sqrt\ga\H^{d-1}(\G)
\end{align*}
and
$$\H^{d-1}\Big(
\G
\smallsetminus\bigcup_{i\in I_1\cup I_2\cup I_3\cup I_4}B(x_i,r_i(1-2\sqrt\ga))
\Big)\,<\,\gamma+
4d\sqrt\ga\H^{d-1}(\G)\,.$$
We have a finite number of disjoint closed balls
$B(x_i,r_i(1-2\sqrt\ga))$, $i\in I_1\cup I_2\cup I_3\cup I_4$.
By increasing slightly all the radii $r_i$,
we can keep the balls disjoint, 
ensure that each radius $r_i$ satisfies
the same strict inequalities
for $i$ in 
$I_1\cup I_2\cup I_3\cup I_4$,
and get the inequality
$$\H^{d-1}\Big(
\G
\smallsetminus\bigcup_{i\in I_1\cup I_2\cup I_3\cup I_4}\bro(x_i,r_i(1-2\sqrt\ga))
\Big)\,<\,2\gamma+4d\sqrt\ga\H^{d-1}(\G)\,.$$
The above set is a compact subset of $\G$.
For $k=1,2$, we define
$$R_k\,=\,\overline{\G}_k
\smallsetminus\bigcup_{i\in I_1\cup I_2\cup I_3\cup I_4}\bro(x_i,r_i(1-2\sqrt\ga))
\,.$$
The sets $R_1$ and $R_2$ 
are compact and their $\H^{d-1}$ measure is less than
$2\gamma+4d\sqrt\ga\H^{d-1}(\G)$
(recall that $\partial_\G\G^1$ and
$\partial_\G\G^2$ have a null $\H^{d-1}$ measure).
For $k=1,2$,
by the definition of the Hausdorff measure~$\H^{d-1}$, there exists
a
collection of balls
$B(y_j,s_j)$, $j\in J_k$ 
such that:
$$\displaylines{
\forall j\in J_k\qquad 
0<s_j<\min\big(\eta_0,\frac{r_\G}{2}
\big)
\,,\qquad
B(y_j,s_j)\cap R_k\neq\varnothing\,,\cr
\sum_{j\in J_k}\al_{d-1} s_j^{d-1}\,<\,3\gamma+4d\sqrt\ga\H^{d-1}(\G)
\,,\cr
R_k
\,\subset\,
\bigcup_{j\in J_k}\bro(y_j,s_j)\,.
}$$
By compactness of $R_1$ and $R_2$, 
the sets $J_1$ and $J_2$ can be chosen to be finite.
It remains to cover
$$R_0\,=\,\G
\smallsetminus\bigcup_{i\in I_1\cup I_2\cup I_3\cup I_4}\bro(x_i,r_i(1-2\sqrt\ga))
\smallsetminus
\bigcup_{j\in J_1\cup J_2}\bro(y_j,s_j)\,.$$
The set $R_0$
is a closed subset of $\G$
which is at a positive distance from $\G^1$
and $\G^2$.
There exists
a
collection of balls
$B(y_j,s_j)$, $j\in J_0$ 
such that:
$$\displaylines{
\forall j\in J_0\qquad 
0<s_j<\min\big(\eta_0,\frac{r_\G}{2},
\frac{1}{8d}\,\text{dist}(R_0,
\G^1\cup\G^2)
\big)
\,,\qquad
B(y_j,s_j)\cap R_0\neq\varnothing\,,\cr
\sum_{j\in J_0}\al_{d-1} s_j^{d-1}\,<\,3\gamma+4d\sqrt\ga\H^{d-1}(\G)
\,,\cr
R_0\,\subset\,
\bigcup_{j\in J_0}\bro(y_j,s_j)\,.
}$$
Now the collection of balls
$$
\bro(x_i,r_i(1-2\sqrt\ga)),\,
{i\in I_1\cup I_2\cup I_3\cup I_4},\quad
B(y_j,s_j), \,j\in J_0\cup J_1\cup J_2$$ 
covers completely $\G$. We will next replace these balls by
polyhedra.
For $j\in J_0\cup J_1\cup J_2$, let $x_j$ belong to
$B(y_j,s_j)\cap \G$ and let $Q_j$ be the cube $Q(x_j,4s_j)$.
For $i$ in $I_1\cup I_2\cup I_3\cup I_4$, the point $x_i$
belongs to exactly one hypersurface among $S_1,\dots, S_p$,
which we denote by
$S_{s(i)}$. In particular
$\G$ admits a normal vector $v_\O(x_i)$
at $x_i$ in the classical sense. 
For each~$i$ in~$I_1\cup I_2\cup I_3\cup I_4$, 
let $P_i$ be a convex open polygon inside the 
hyperplane $\hyp(x_i,v_\O(x_i))$ such that
$$\displaylines{\disc(x_i,r_i(1-2\sqrt\ga),v_\O(x_i))\subset P_i
\subset \disc(x_i,r_i(1-\sqrt\ga),v_\O(x_i))\,,\cr
|\H^{d-2}(\partial P_i)-\al_{d-2} r_i^{d-2}
(1-\sqrt\ga)^{d-2}|\,\leq\,\delta_0\al_{d-2} r_i^{d-2}(1-\sqrt\ga)^{d-2}
\,,\quad\cr
|\H^{d-1}(P_i)-\al_{d-1} r_i^{d-1}
(1-\sqrt\ga)^{d-1}|\leq \delta_0 \al_{d-1} r_i^{d-1}
(1-\sqrt\ga)^{d-1}
\,.}$$
Thanks to the choices
of the radius $r_i$ and the constants $M_0,\eta_0$, we have then
$$\displaylines{
\G\cap B(x_i,r_i(1-2\sqrt\ga))\,
\subset\, S_{s(i)}\cap B(x_i,r_i(1-2\sqrt\ga))\,\subset\, 
\smash{
\cyl
}^{\kern-5pt\!\!\!\raise5pt\hbox{\scriptsize o}}
\kern5pt
(P_i,2\ga r_i) 
\,,\cr
\G\cap B(x_i,r_i)\,
\subset\, S_{s(i)}\cap B(x_i,r_i)\,\subset\, 
\cyl( \disc(x_i,r_i,v_\O(x_i)) ,M_0\delta_0 r_i)
\,,\cr
\forall x\in B(x_i,r_i)\cap\G\qquad
|v_\O(x)-v_\O(x_i)|_2\,<\,1\,.}$$
The choice of $\delta_0$ guarantees that
$M_0\delta_0(1+\delta_0)r_i<2\ga r_i$.
Let $t$ be such that
$$M_0\delta_0(1+\delta_0)r_i\leq t<\sqrt\ga r_i\,.$$
We have
$$\displaylines{
-tv_{\O}(x_i)+P_i\subset \O\cap B(x_i,r_i)\,,\qquad
\G\cap 
(-tv_{\O}(x_i)+P_i)=\varnothing\,.}$$
In particular, the set $\G$ can intersect the cylinder
$\cyl(P_i,t)$ 
only along its
lateral sides, which are parallel to $v_\O(x_i)$. Let 
$x$ belong to 
$\G\cap\partial \cyl(P_i,t)$. Then  
$$|v_{\cyl(P_i,t)}(x)- v_{\O}(x)|_2\,\geq\,
|v_{\cyl(P_i,t)}(x)- v_{\O}(x_i)|_2-
|v_{\O}(x_i)- v_{\O}(x)|_2\,\geq\,\sqrt{2}-1\,.$$
Therefore the cylinder 
$\cyl(P_i,t)$
is transverse to~$\G$.
We will replace the ball $\bro(x_i,r_i(1-2\sqrt\ga))$
by the cylinder
$\cyl(P_i,t_i)$, for a carefully chosen value
of $t_i$ in the interval
$[M_0\delta_0(1+\delta_0)r_i,\sqrt\ga r_i[$.
However, we must
delay the choices of the values $t_i$,
$i\in I_3\cup I_4$ 
until we have modified the set $F$ inside $\O$.
We deal next with the interfaces inside $\O$ and we make 
an approximation of $F$
controlled by a factor $\eps$. We choose $\eps$ sufficiently small
compared to $\ga$ so that, when we perturb the set $F$ by 
a volume $\eps$, the resulting effect close to $\G$
is still of order $\ga$. 
Let $\eps$ be such that $0<\eps<\gamma$ and
$$\eps \,<\,\ga\al_d\min_{i\in I_1\cup I_2\cup I_3\cup I_4}r_i^d\,.$$
We use next a classical approximation result:
there exists a relatively closed subset $L$ of $\O$ having finite
perimeter such that $\O\cap\partial L$ is an hypersurface
of class $C^\infty$ and
$$\ld(F\Delta L)\,<\,\eps\,,\qquad
\Big|\int_{\O\cap\rb F}
\kern-20pt
\nu(v_{F}(y))\,d\H^{d-1}(y)
-
\int_{\O\cap\partial L}
\kern-20pt
\nu(v_{L}(y))\,d\H^{d-1}(y)\Big|\,<\,\eps\,.$$
In the case where $\nu$ is constant, this result is stated
in Lemma 4.4 of \cite{QG}. In the non constant case, the
argument should be slightly modified, as explained in the proof
of proposition 14.8 of \cite{Cerf:StFlour}, where the approximation
is performed in $\RR^d$ instead of $\O$. When working inside $\O$,
the extra difficulty is to deal with regions close
to the boundary (see the proof of Proposition 4.3 of \cite{QG}).
For $r>0$, we define
$$\partial L_r\,=\,\big\{\,x\in\partial L:d(x,\G)\geq r\,\big\}\,.$$
By continuity of the
measure $\H^{d-1}|_{\partial L}$,
there exists $r^*>0$ such that
$$\H^{d-1}(\O\cap\partial L\smallsetminus \partial L_{2r^*})
\,\leq\,
 \eps\,.$$
We apply lemma~\ref{chapitre7rfpres} to the set 
$\partial L_{r^*}$
and the hypersurface $\O\cap\partial L$:
$$\displaylines{
\exists M>0\quad
\forall\delta>0\quad\exists\,\eta>0\quad
\forall x,y\in 
\partial L_{r^*}
\quad
|x-y|_2\leq \eta\,\Rightarrow\,
d_2\big(y,\tan(\partial L,x)\big)
\leq M\delta |x-y|_2\,.}$$
For a point~$x$ belonging to 
$\partial L_{r^*}$,
the tangent hyperplane of~$\O\cap\partial L$ at~$x$ 
is precisely $\hyp(x,v_L(x))$.
Let $M$ be as above. We can assume that $M>1$.
Let $\delta$ 
in $]0,\delta_0[$ be such that $2\delta M<\eps$.
Let $\eta$ be associated to $\delta$ as in the above property.
For $x\in \partial L_{2r^*}$,
$$\displaylines{
\lim_{r\to 0}\quad
(\al_{d-1}r^{d-1})^{-1}
\H^{d-1}(B(x,r)\cap\partial L)\,=\,1\,,\cr
\lim_{r\to 0}\quad
(\al_{d-1}r^{d-1})^{-1}
\int_{B(x,r)\cap\partial L}
\kern-20pt
\nu(v_L(y))\,d\H^{d-1}(y)
\,=\,\nu(v_L(x))\,.}$$
For any~$x$ in~
$\partial L_{2r^*}$,
there exists a positive~$r(x,\eps)$ such that, 
for any $r<r(x,\eps)$, 
$$\displaylines{
|\H^{d-1}(B(x,r)\cap \partial L)-\al_{d-1}r^{d-1}|
\,\leq\,\eps\,\al_{d-1}r^{d-1}
\,,\cr
\Big|
(\al_{d-1}r^{d-1})^{-1}
\int_{B(x,r)\cap \partial L}
\kern-20pt
\nu(v_L(y))\,d\H^{d-1}(y)
-\nu(v_L(x))
\Big|\,\leq\,
\eps\,
\,.}$$
The family of balls
$B(x,r)$,
$x\in 
\partial L_{2r^*}$,
$r<\min(r^*,\eta_0,
r(x,\eps),\eps,\eta)$,
is a Vitali relation for 
$\partial L_{2r^*}$.
By the standard Vitali covering Theorem,
we may select a finite or countable collection
of disjoint balls $B(x_i,r_i)$, $i\in I'$, such that:
for any~$i$ in~$I'$, 
$x_i\in \partial L_{2r^*}$,
$$r_i\,<\,\min(r^*,\eta_0,r(x_i,\eps),\eps,\eta)$$
and
$$\text{either}\qquad\H^{d-1}\Big(
\partial L_{2r^*}
\smallsetminus\bigcup_{i\in I'}B(x_i,r_i)
\Big)\,=\,0\qquad\text{or}\qquad
\sum_{i\in I'}r_i^{d-1}\,=\,\infty
\,.$$
Because for each~$i$ in~$I'$, $r_i$ is smaller than 
$r(x_i,\eps)$, 
$$\al_{d-1}(1-\eps)\sum_{i\in I'}r_i^{d-1}\,
\leq\,\H^{d-1}(\O\cap\partial L)\,<\,\infty\,$$
and therefore the first case occurs, so that we may select a finite
subset~$I_5$ of~$I'$ such that
 $$\H^{d-1}\Big(
\partial L_{2r^*}
\smallsetminus\bigcup_{i\in I_5}B(x_i,r_i)
\Big)\,<\,\eps\,.$$
We have a finite number of disjoint closed balls
$B(x_i,r_i)$, $i\in I_5$. By increasing slightly all the radii $r_i$,
we can keep the balls disjoint, each $r_i$ strictly smaller than
$\min(r^*,\eta_0,r(x_i,\eps),\eps,\eta)$ for $i$ in $I_5$,
and get the stronger inequality
$$\H^{d-1}\Big(
\partial L_{2r^*}
\smallsetminus\bigcup_{i\in I_5}\bro(x_i,r_i)
\Big)\,<\,\eps\,.$$
For each~$i$ in~$I_5$, let $P_i$ be a convex open polygon inside the 
hyperplane $\hyp(x_i,v_L(x_i))$ such that
$$\displaylines{\disc(x_i,r_i,v_L(x_i))\subset P_i
\subset \disc(x_i,r_i(1+\delta),v_L(x_i))\,,\cr
|\H^{d-2}(\partial P_i)-\al_{d-2} r_i^{d-2}|\,\leq\,\delta\al_{d-2} r_i^{d-2}
\,,\cr
|\H^{d-1}(P_i)-\al_{d-1} r_i^{d-1}|\leq \delta \al_{d-1} r_i^{d-1}\,.}$$
We set $\psi=M\delta(1+\delta)$ (hence $\psi<\eps<1$).
Let $i$ belong to~$I_5$.
Let
$D_i$ be the cylinder 
$$D_i\,=\,\cyl(P_i,M\delta(1+\delta) r_i)$$
of basis $P_i$ and height
$2\psi r_i$.
The point $x_i$ belongs to $\partial L_{2r^*}$, 
the radius $r_i$ is smaller than~$\eta$ and $r^*$, so that
$$\displaylines{
\forall x\in\partial L
\cap B(x_i,r_i)\quad
d_2\big(x,\hyp(x_i,v_L(x_i))\big)
\leq M\delta|x-x_i|_2\,,}$$
whence
$$\partial L\cap B(x_i,r_i)\,\subset\,
\cyl\big(\disc(x_i,r_i,v_L(x_i)),M\delta r_i\big)
\,\subset\, \dro_i\,.$$
We will approximate $F$ by $L$ inside $\O$ and we will
push the interfaces $\G^1 \cap \p^*(\O\smallsetminus F)$ and $\G^2 \cap
\p^* F$
into $\O$. 
We next handle the regions close to $\G$ 
inside the family of
balls $B(x_i,r_i)$, $i\in I_1\cup I_2\cup I_3\cup I_4$. 
We will modify adequately
the set $F$ to ensure that no significant 
interface is created within these balls. 
Our technique consists in building a small flat cylinder centered 
on $\G$ which we add (for indices in $I_1\cup I_3$)
or remove (for indices in $I_2\cup I_4$)
to the set $F$.
We have to design carefully this operation in order not to create
any significant additional interface. This is the place where we tie
together the covering of the boundary and the inner approximation.
Recall that we already chose a family of polygons $P_i$, 
$i\in I_1\cup I_2\cup I_3\cup I_4$. 
For 
$i \in I_1\cup I_2$, we simply define $D_i$ to be the cylinder
$$D_i\,=\,
\cyl(P_i,
M_0\delta_0(1+\delta_0)r_i)\,,$$
see figure \ref{chapitre7polygone3}.
The construction of the cylinders associated to the indices
$i\in I_3\cup I_4$ is more complicated.
Our technique consists in choosing carefully the height $t_i$
of the cylinders $\cyl(P_i,t_i)$ for $i\in I_3\cup I_4$.
We examine separately the indices in $I_3$ and $I_4$.

{\noindent {\bf $\bullet$ Balls indexed by $I_3$.}}
Let $i$ belong to $I_3$.
Because of the condition imposed on $\eps$, 
we have
$$\displaylines{
|\ld(B(x_i,r_i)\cap L)-\al_{d}r_i^{d}/2|
\,\leq\,\gamma\,\al_{d}r_i^{d}+\eps
\,\leq\,2\gamma\,\al_{d}r_i^{d}\,.
}$$
Since in addition
$$|\ld(B(x_i,r_i)\smallsetminus\O)-\al_{d}r_i^{d}/2|
\,\leq\,\gamma\,\al_{d}r_i^{d}\,,
$$
it follows that
$$\ld(B(x_i,r_i)\cap(\O\smallsetminus \lro))
\,\leq\,3\gamma\,\al_{d}r_i^{d}\,.
$$
Thanks to the choice of the polygon~$P_i$, we have then
$$\displaylines{
\int_{2\ga r_i<t<\sqrt\ga r_i}
\H^{d-1}( (-tv_{\O}(x_i)+P_i)\smallsetminus \lro)
\,dt
\,\leq\,
\ld(B(x_i,r_i)\cap(\O\smallsetminus \lro))
\,\leq\,
3\ga\al_d r_i^d\,.}$$
The condition on $\ga$ yields in particular
$\sqrt\ga-2\ga\geq\sqrt\ga/2$.
Hence there exists $t_i\in ]2\ga r_i,\sqrt\ga r_i[$
such that
$$
\H^{d-1}( (-t_iv_{\O}(x_i)+P_i)\smallsetminus \lro)
\,\leq\,
6\sqrt\ga\al_d r_i^{d-1}\,.$$
Let $D_i$ be the cylinder $D_i=\cyl(P_i,t_i)$.

\begin{figure}
\centering
\begin{picture}(0,0)%
\includegraphics{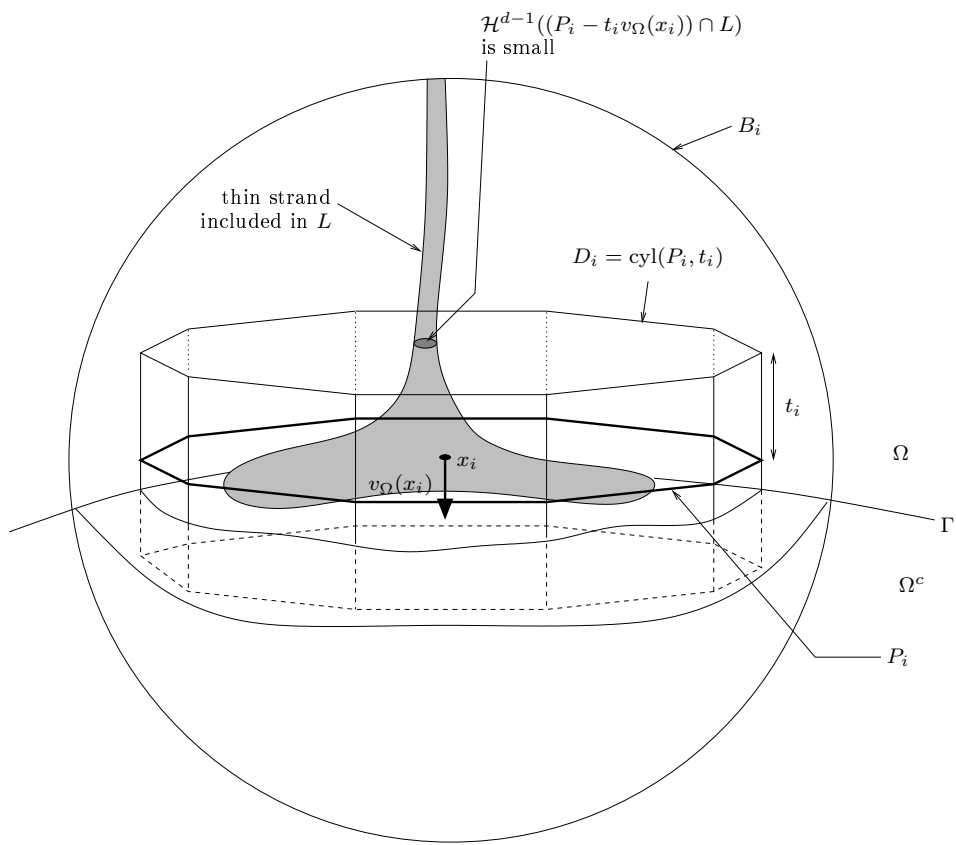}%
\end{picture}%
\setlength{\unitlength}{1973sp}%
\begingroup\makeatletter\ifx\SetFigFont\undefined%
\gdef\SetFigFont#1#2#3#4#5{%
  \reset@font\fontsize{#1}{#2pt}%
  \fontfamily{#3}\fontseries{#4}\fontshape{#5}%
  \selectfont}%
\fi\endgroup%
\begin{picture}(11727,10438)(439,-9968)
\put(6076,-5236){\makebox(0,0)[lb]{\smash{{\SetFigFont{8}{9.6}{\rmdefault}{\mddefault}{\updefault}{\color[rgb]{0,0,0}$x_i$}%
}}}}
\put(11551,-5161){\makebox(0,0)[lb]{\smash{{\SetFigFont{8}{9.6}{\rmdefault}{\mddefault}{\updefault}{\color[rgb]{0,0,0}$\O$}%
}}}}
\put(11626,-6811){\makebox(0,0)[lb]{\smash{{\SetFigFont{8}{9.6}{\rmdefault}{\mddefault}{\updefault}{\color[rgb]{0,0,0}$\O^c$}%
}}}}
\put(12151,-6061){\makebox(0,0)[lb]{\smash{{\SetFigFont{8}{9.6}{\rmdefault}{\mddefault}{\updefault}{\color[rgb]{0,0,0}$\G$}%
}}}}
\put(9601,-1036){\makebox(0,0)[lb]{\smash{{\SetFigFont{8}{9.6}{\rmdefault}{\mddefault}{\updefault}{\color[rgb]{0,0,0}$B_i$}%
}}}}
\put(11476,-7711){\makebox(0,0)[lb]{\smash{{\SetFigFont{8}{9.6}{\rmdefault}{\mddefault}{\updefault}{\color[rgb]{0,0,0}$P_i$}%
}}}}
\put(10201,-4561){\makebox(0,0)[lb]{\smash{{\SetFigFont{8}{9.6}{\rmdefault}{\mddefault}{\updefault}{\color[rgb]{0,0,0}$t_i$}%
}}}}
\put(5776,-5536){\makebox(0,0)[rb]{\smash{{\SetFigFont{8}{9.6}{\rmdefault}{\mddefault}{\updefault}{\color[rgb]{0,0,0}$v_{\O}(x_i)$}%
}}}}
\put(4501,-1936){\makebox(0,0)[rb]{\smash{{\SetFigFont{8}{9.6}{\rmdefault}{\mddefault}{\updefault}{\color[rgb]{0,0,0}thin strand}%
}}}}
\put(4501,-2236){\makebox(0,0)[rb]{\smash{{\SetFigFont{8}{9.6}{\rmdefault}{\mddefault}{\updefault}{\color[rgb]{0,0,0}included in $L$}%
}}}}
\put(6376,239){\makebox(0,0)[lb]{\smash{{\SetFigFont{8}{9.6}{\rmdefault}{\mddefault}{\updefault}{\color[rgb]{0,0,0}$\H^{d-1} ((P_i -t_i v_{\O}(x_i)) \cap L )$}%
}}}}
\put(6376,-61){\makebox(0,0)[lb]{\smash{{\SetFigFont{8}{9.6}{\rmdefault}{\mddefault}{\updefault}{\color[rgb]{0,0,0}is small}%
}}}}
\put(8476,-2686){\makebox(0,0)[b]{\smash{{\SetFigFont{8}{9.6}{\rmdefault}{\mddefault}{\updefault}{\color[rgb]{0,0,0}$D_i = \cyl(P_i, t_i)$}%
}}}}
\end{picture}%
\caption{The cylinder $D_i$ for $i\in I_4$.}
\label{chapitre7polygone4}
\end{figure}

{\noindent {\bf $\bullet$ Balls indexed by $I_4$.}}
Let $i$ belong to $I_4$.
Because of the condition imposed on $\eps$, 
we have
$$\displaylines{
|\ld(B(x_i,r_i)\cap(\O\smallsetminus L))-\al_{d}r_i^{d}/2|
\,\leq\,\gamma\,\al_{d}r_i^{d}+\eps
\,\leq\,2\gamma\,\al_{d}r_i^{d}\,.
}$$
Since in addition
$$|\ld(B(x_i,r_i)\smallsetminus\O)-\al_{d}r_i^{d}/2|
\,\leq\,\gamma\,\al_{d}r_i^{d}\,,
$$
it follows that
$$\ld(B(x_i,r_i)\cap L)
\,\leq\,3\gamma\,\al_{d}r_i^{d}\,.
$$
Thanks to the choice of the polygon~$P_i$, we have then
$$\displaylines{
\int_{2\ga r_i<t<\sqrt\ga r_i}
\H^{d-1}( (-tv_{\O}(x_i)+P_i)\cap L)
\,dt
\,\leq\,
\ld(B(x_i,r_i)\cap L)
\,\leq\,
3\ga\al_d r_i^d\,.}$$
The condition on $\ga$ yields in particular
$\sqrt\ga-2\ga\geq\sqrt\ga/2$.
Hence there exists $t_i\in ]2\ga r_i,\sqrt\ga r_i[$
such that
$$
\H^{d-1}( (-t_iv_{\O}(x_i)+P_i)\cap L)
\,\leq\,
6\sqrt\ga\al_d r_i^{d-1}\,.$$
Let $D_i$ be the cylinder $D_i=\cyl(P_i,t_i)$ (see figure
\ref{chapitre7polygone4}).
We have now built the whole family of cylinders
$D_i$,
${i\in I_1\cup I_2\cup I_3\cup I_4\cup I_5}$.
Moreover, the sets
$$\dro_i\,,\quad i\in I_1\cup I_2\cup I_3\cup I_4\,,\qquad
\bro(y_j,s_j)\,,\quad j\in J_0\cup J_1\cup J_2\,,$$
cover completely $\G$.
It remains now to cover the region
$$R_3\,=\,\O\cap\partial L\kern 3pt\smallsetminus
\kern-7pt
\bigcup_{i\in I_1\cup I_2\cup I_3\cup I_4\cup I_5}
\kern-11pt
\dro_i
\kern 5pt
\kern-3pt
\smallsetminus\bigcup_{j\in J_0\cup J_1\cup J_2}
\kern-7pt
\bro(y_j,s_j)\,.
$$
Since $R_3$ does not intersect~$\G$, the distance
$$\rho=\frac{1}{8d}\,\text{dist}
(
\G,R_3
)$$
is positive and also $R_3$ is compact.
From the preceding inequalities, we deduce that
$$\displaylines{\H^{d-1}(
R_3
)
\,\leq\,
\H^{d-1}(\O\cap\partial L\smallsetminus \partial L_{2r^*})+
\H^{d-1}\Big(\partial L_{2r^*}\smallsetminus\bigcup_{i\in I_5}\dro_i\Big)
\cr
\,\leq\,
\eps+
\H^{d-1}\Big(\partial L_{2r^*}\smallsetminus\bigcup_{i\in I_5}\bro(x_i,r_i)\Big)
\,\leq\,
2\eps\,.
}$$
By the definition of the Hausdorff measure~$\H^{d-1}$, there exists
a collection of balls
$B(y_j,s_j)$, $j\in J_3$, such that:
$$\displaylines{
\forall j\in J_3\qquad 0<s_j<\rho,\qquad
B(y_j,s_j)\cap
R_3
\neq\varnothing\,,\cr
R_3
\,\subset\,
\bigcup_{j\in J_3}\bro(y_j,s_j)\,,\cr
\sum_{j\in J_3}\al_{d-1} s_j^{d-1}\,\leq\,
3\eps\,.
}$$
By compactness, we might assume in addition that
$J_3$ is finite.
For $j\in J_3$, let $x_j$ belong to
$B(y_j,s_j)\cap
R_3
$
and let $Q_j$ be the cube $Q(x_j,4s_j)$.
We set 
$$\displaylines{P\,=\,\bigg((\O\cap L)
\cup
\bigcup_{i\in I_1\cup I_3\cup I_5} D_i
\cup \bigcup_{j\in J_1} Q_j
\bigg)
\smallsetminus
\bigcup_{i\in I_2\cup I_4} D_i
\smallsetminus
\bigcup_{j\in J_0\cup J_2\cup J_3} Q_j\,.
}$$
The sets 
$\qro_j$, $j\in J_0\cup J_1\cup J_2\cup J_3$, 
$\dro_i$, $i\in I_1\cup I_2\cup I_3\cup I_4\cup I_5$
cover $\partial L\cup\G$, therefore
$$\displaylines{
\partial P
\,\,\subset\,\,
\bigcup_{i\in I_1\cup I_2\cup I_3\cup I_4\cup I_5}
\partial D_i
\cup\bigcup_{j\in J_0\cup J_1\cup J_2\cup J_3}\partial Q_j\,,
}$$
thus $P$ is polyhedral and $\partial P$ is 
transverse to $\G$.
Since the sets
$$\dro_i\,,\quad i\in I_1\cup I_3\,,\qquad
\qro_j\,, \quad j\in J_1$$
cover completely $\overline \G^1$, 
while the sets
$$D_i\,,\quad i\in I_2\cup I_4\cup I_5\,,\qquad
Q_j\,, \quad j\in J_0\cup J_2\cup J_3$$
do not intersect
$\overline \G^1$, 
then
$\overline \G^1$ is included in the interior of~$P$.
Similarly, the sets
$$\dro_i\,,\quad i\in I_2\cup I_4\,,\qquad
\qro_j\,, \quad j\in J_2$$
cover completely $\overline \G^2$, 
while the sets
$$D_i\,,\quad i\in I_1\cup I_3\cup I_5\,,\qquad
Q_j\,, \quad j\in J_0\cup J_1\cup J_3$$
do not intersect
$\overline \G^2$, 
thus
$\overline \G^2$ is included in the interior of the complement
of~$P$.
We next check that the set $P\cap \O$
approximates the initial set $F$
with respect to the volume. We have
$$(P \cap \O) \Delta F\,\subset\,
(L\Delta F)
\cup
\bigcup_{i\in I_1\cup I_2\cup I_3\cup I_4\cup I_5}D_i
\cup
\bigcup_{j\in J_0\cup J_1\cup J_2\cup J_3}Q_j$$
whence
$$\displaylines{
\ld((P\cap \O)\Delta F)\,\leq\,
\eps+
\hfill\cr
\sum_{i\in I_1\cup I_2\cup I_3\cup I_4}
2\al_{d-1} r_i^{d-1}(1+\delta_0)\sqrt\ga r_i +
\sum_{i\in I_5}2\al_{d-1} r_i^{d-1} (1+\delta)\psi r_i
+\sum_{j\in J_0\cup J_1\cup J_2\cup J_3} \al_d (2s_j)^d\,.}$$
Yet each $r_i$ is smaller than $\ga$,
$$\displaylines{
\sum_{i\in I_1\cup I_2\cup I_3\cup I_4}\al_{d-1} r_i^{d-1}\leq 2\H^{d-1}(\G)
\,,\cr
\sum_{i\in I_5}\al_{d-1} r_i^{d-1}\leq 2\H^{d-1}(\O\cap\partial L)
\leq \frac{2}{\nu_{\text{min}}}(\nu_{\text{max}}\H^{d-1}(\p^* F \cap \O)+\eps)\,,
\cr
\sum_{j\in J_0\cup J_1\cup J_2\cup J_3}\al_{d-1} s_j^{d-1}\leq 
3\big(3\gamma+4d\sqrt\ga\H^{d-1}(\G)\big)+3\eps\,,}$$
so that
\begin{align*}
\ld((P\cap \O)\Delta F) & \,\leq\,
\eps
+6\sqrt\gamma\H^{d-1}(\G) +
\frac{6\eps}{\nu_{\text{min}}}(\nu_{\text{max}}\H^{d-1}(\p^* F\cap \O)+\eps)\\
& \qquad \qquad +3\cdot 2^d
\frac{\al_d}{\al_{d-1}}(3\gamma+4d\sqrt\ga\H^{d-1}(\G)+\eps)\,.
\end{align*}
We 
estimate next 
the capacity
of~$P$.
To do this, we examine the intersection of $\p P \cap \O$ with each
polyhedral cylinder. For $i\in I_1\cup I_2$, we use the obvious inclusion
$$\ P \cap \O \cap\partial D_i\,\subset\,
\O\cap\partial D_i\,.$$
For $i\in I_3\cup I_4\cup I_5$, 
the sets
$\p P \cap \O\cap\partial D_i$
require more attention.
We consider separately the indices of $I_3$, $I_4$ and $I_5$.

{\noindent {\bf $\bullet$ Cylinders indexed by $I_3$.}}
Let 
$i$ in~$I_3$. We have
$$\O\cap\partial P\cap \partial D_i\,\subset\,
\O\cap(\partial D_i\smallsetminus\lro)
\cup\bigcup_{j\in J_0\cup J_1\cup J_2\cup J_3}\partial Q_j\,.$$
Yet, thanks to the construction of the cylinder $D_i$,
$$\displaylines{\H^{d-1}(\O\cap\partial D_i\smallsetminus \lro)\,\leq\,
\H^{d-1}( (-t_iv_{\O}(x_i)+P_i)\smallsetminus \lro)+
\H^{d-2}(\partial P_i)2\sqrt\ga r_i
\hfill\cr
\,\leq\,
6\sqrt\ga\al_d r_i^{d-1}
+2\al_{d-2} r_i^{d-2}2\sqrt\ga r_i
\,\leq\,
6\sqrt\ga(\al_d+
\al_{d-2})
r_i^{d-1}
\,.}$$

{\noindent {\bf $\bullet$ Cylinders indexed by $I_4$.}}
Let 
$i$ in~$I_4$. We have
$$\O\cap\partial P\cap \partial D_i\,\subset\,
\O\cap(\partial D_i\cap L)
\cup\bigcup_{j\in J_0\cup J_1\cup J_2\cup J_3}\partial Q_j\,.$$
Yet, thanks to the construction of the cylinder $D_i$,
$$\displaylines{\H^{d-1}(\O\cap\partial D_i\cap L)\,\leq\,
\H^{d-1}( (-t_iv_{\O}(x_i)+P_i)\cap L)+
\H^{d-2}(\partial P_i)2\sqrt\ga r_i
\hfill\cr
\,\leq\,
6\sqrt\ga\al_d r_i^{d-1}
+2\al_{d-2} r_i^{d-2}2\sqrt\ga r_i
\,\leq\,
6\sqrt\ga(\al_d+
\al_{d-2})
r_i^{d-1}
\,.}$$

{\noindent {\bf $\bullet$ Cylinders indexed by $I_5$.}}
Let 
$i$ in~$I_5$.
We set
$$G_i\,=\,
\disc\big(
x_i-\psi r_i v_{L}(x_i),
\sqrt{1-\psi^2} 
 r_i,v_L(x_i)\big)
\,.$$
We claim that the set $G_i$
is included in the interior of~$L$.
Indeed, 
$G_i\subset
B(x_i,r_i)\cap \partial D_i$, yet
$\partial L\cap B(x_i,r_i)\subset\dro_i$,
therefore
$G_i$ does not intersect~$\partial L$. 
Since $v_{L}(x_i)$ is 
the exterior normal vector to~$L$ 
at~$x_i$,
then $G_i$ is included in 
$\lro$.
The definition of the set $P$ implies that
$$\partial P\cap G_i
\,\subset\,\bigcup_{j\in J_0\cup J_1\cup J_2\cup J_3}\partial Q_j\,,$$
whence
$$\O\cap\partial P\cap\partial D_i\,\subset\,
(\partial D_i\smallsetminus G_i)\cup
\,\bigcup_{j\in J_0\cup J_1\cup J_2\cup J_3}\partial Q_j\,.$$
Yet
$$\displaylines{
\H^{d-1}\big(\partial D_i
\smallsetminus
(P_i+\psi r_i v_{L}(x_i))
\smallsetminus G_i
\big)\,\leq\,
2\al_{d-2} r_i^{d-2}2\psi r_i+
\al_{d-1} r_i^{d-1}
\big(1+\delta-
(1-\psi^2)^{(d-1)/2}\big)
\cr
\,\leq\,
\al_{d-1} r_i^{d-1}
\Big(
4\frac{\al_{d-2}}{\al_{d-1}} \psi +
1+\delta-
(1-\psi^2)^{(d-1)/2}\Big)
\,.
}$$
Finally, we conclude that
\begin{align*}
\O\cap\partial P
 & \,\,\subset\,\,
\bigcup_{i\in I_1\cup I_2}
(\O\cap\partial D_i)
\cup
\bigcup_{i\in I_3}
(\O\cap D_i\smallsetminus\lro)
\cup
\bigcup_{i\in I_4}
(\O\cap\partial D_i\cap L) \\
& \qquad \qquad \cup\bigcup_{i\in I_5}
(\partial D_i\smallsetminus G_i)
\cup
\kern -7pt
\bigcup_{j\in J_0\cup J_1\cup J_2\cup J_3}
\kern -7pt
\partial Q_j\,.
\end{align*}
Therefore
\begin{align*}
\I_{\O}(P) & \,\leq\,
\sum_{i\in I_1\cup I_2}
\int_{\O\cap\partial D_i}
\kern-20pt
\nu(v_P(x))\,d\H^{d-1}(x)
+\nu_{\max}
\sum_{i\in I_3}
\H^{d-1}(\O\cap\partial D_i\smallsetminus \lro)\\
& \quad +\nu_{\max}
\sum_{i\in I_4}
\H^{d-1}(\O\cap\partial D_i\cap L)
\hfill\cr \\
 & \quad +\sum_{i\in I_5}
\Big(
\nu(v_L(x_i))\H^{d-1}(P_i)+
\nu_{\max}
\H^{d-1}\big(\partial D_i
\smallsetminus
(P_i+\psi r_i v_{L}(x_i))
\smallsetminus G_i
\big)
\Big)\\
& \quad +\nu_{\max}
\sum_{j\in J_0\cup J_1\cup J_2\cup J_3}\H^{d-1}(\partial Q_j)\,. 
\end{align*}
We use now the various estimates obtained in the course
of the approximation. 
We get
\begin{align*}
 & \I_{\O}(P) \, \leq\,
\sum_{i\in I_1\cup I_2}
\Big(\al_{d-1} r_i^{d-1}(1+\delta_0)\nu(v_\O(x_i))+
\nu_{\max} \al_{d-2} r_i^{d-1}2M_0\de_0(1+\de_0)^2\Big)\\
& \quad +
\sum_{i\in I_3\cup I_4}
\nu_{\max}\Big(
6\sqrt\ga(\al_d+
\al_{d-2})
r_i^{d-1}
\Big) \\
& \quad +
\sum_{i\in I_5}
\Big(\al_{d-1} r_i^{d-1}(1+\delta)\nu(v_L(x_i))\\
& \qquad+
\nu_{\max} 
\al_{d-1} r_i^{d-1}
\Big(
4\frac{\al_{d-2}}{\al_{d-1}} \psi +
1+\delta-
(1-\psi^2)^{(d-1)/2}\Big)
\Big) \\
& \quad +\sum_{j\in  J_0\cup J_1\cup J_2\cup J_3}
\nu_{\max}
\al_{d-1}2^{d-1} s_j^{d-1}
\hfill\\
& \leq\,
\frac{1+\de_0}{1-\ga}
\sum_{i\in I_1}
\int_{B(x_i,r_i)\cap\rb(\O\smallsetminus F)}\nu({v_\O}(y))\,d\H^{d-1}(y)\\
&\quad+ \frac{1+\de_0}{1-\ga}
\sum_{i\in I_2}
\int_{B(x_i,r_i)\cap\rb F}\nu({v_\O}(y))\,d\H^{d-1}(y)
\\
&\quad +
\frac{1+\de}{ 1-\eps}
\sum_{i\in I_5}
\int_{B(x_i,r_i)\cap\partial L}\nu({v_L}(y))\,d\H^{d-1}(y)\\
& \quad+\sum_{i\in I_1\cup I_2\cup I_3\cup I_4\cup I_5}
\nu_{\max} \al_{d-1} r_i^{d-1}
\Big(
\frac{\al_{d-2}}{al_{d-1}}5\ga+
6\sqrt\ga \frac{\al_d+
\al_{d-2}}{\al_{d-1}}
+
4\frac{\al_{d-2}}{\al_{d-1}} \psi \\
& \qquad+ 1+\delta-
(1-\psi^2)^{(d-1)/2}\Big) +\nu_{\max}
2^{d-1}3\big(3\gamma+4d\sqrt\ga\H^{d-1}(\G)+\eps\big)\\
& \leq\,
\frac{1+\de_0}{ 1-\ga}
\bigg(
\int_{\G^1\cap\rb(\O\smallsetminus F)}\nu({v_\O}(y))\,d\H^{d-1}(y)
+
\int_{\G^2\cap\rb F}\nu({v_\O}(y))\,d\H^{d-1}(y)
\\
& \quad +
\int_{\O\cap\partial L}\nu({v_L}(y))\,d\H^{d-1}(y)
\bigg)\\
& \quad +
2(\H^{d-1}(\G)+\H^{d-1}(\O\cap\partial L))
\nu_{\max} 
\Big(
\frac{\al_{d-2}}{\al_{d-1}}5\ga+
6\sqrt\ga\frac{\al_d+
\al_{d-2}}{\al_{d-1}}
+4\frac{\al_{d-2}}{\al_{d-1}} \psi \\
& \qquad +
1+\delta-
(1-\psi^2)^{(d-1)/2}
\Big) +\nu_{\max}\Big(
2^{d-1}3\big(3\gamma+4d\sqrt\ga\H^{d-1}(\G)\big)+3\eps\Big)\\
&\leq\,
\frac{1+\de_0}{ 1-\ga}
\big(\I_{\O}(F)+\eps)\\
& \quad+2\Big(\H^{d-1}(\G)+
\frac{\nu_{\text{max}}
\I_{\O}(F)+\eps}{ \nu_{\text{min}}}
\Big)
\nu_{\max} 
\Big(
\frac{\al_{d-2}}{\al_{d-1}}5\ga+
6\sqrt\ga \frac{\al_d+
\al_{d-2}}{\al_{d-1}}
+\delta\eps+4
\frac{\al_{d-2}}{\al_{d-1}}
\eps\Big)\\
&\quad+\nu_{\max}\Big(
2^{d-1}3\big(3\gamma+4d\sqrt\ga\H^{d-1}(\G)\big)+3\eps\Big) 
\end{align*}
where we have used the inequality $\psi<\eps$ in the last step.
We have also use the inclusions
$$\displaylines{
\forall i\in I_1\qquad B(x_i,r_i)\cap\rb(\O\smallsetminus F)\,\subset\,
\G^1\cap
\rb(\O\smallsetminus F)\,,\cr
\forall i\in I_2\qquad B(x_i,r_i)\cap\rb F\,\subset\,
\G^2\cap
\rb F\,.}$$
Since $\delta_0, \delta,\ga,\eps$ can be chosen arbitrarily small,
we have obtained the desired approximation.\qed


\nocite{ASQU}
\nocite{BE}
\nocite{Ce}
\nocite{DGb}
\nocite{DGCOPI}
\nocite{GI}
\nocite{MAMI}
\nocite{Mattila}


\section{Positivity of $\widetilde{\phi_{\O}}$}
\label{chapitre7secpositif}

We suppose that
\begin{equation}
\label{chapitre7mt1}
 \int_{[0, +\infty[} x \,d\Lambda (x) \,<\,\infty \,,
\end{equation}
We will prove that $\widetilde{\phi_{\O}} >0$ if and only if $\Lambda
(0) < 1-p_c(d)$. In
fact we know that if the condition (\ref{chapitre7mt1}) is satisfied,
$$ \Lambda(0)<1-p_c(d) \quad \Longleftrightarrow \quad \exists v \,,\,\,
\nu(v)>0 \quad \Longleftrightarrow \quad \forall v\,,\,\, \nu(v)>0\,. $$
Thus, the implication
$$ \Lambda(0)\geq 1-p_c(d)\quad \Longrightarrow \quad
\widetilde{\phi_{\O}} =0 $$
is trivial. We suppose that $\Lambda (0) < 1-p_c(d)$. Since $\nu$ satisfies the weak triangle inequality, the function
$v\mapsto \nu(v)$ is continuous, and so as soon as $\Lambda(0)<1-p_c(d)$ and
(\ref{chapitre7mt1}) is satisfied, we have
$$ \nu_{\min} \,=\, \min_{\mathbb{S}^1} \nu \,>\, 0\,.  $$
If $P$ is a polyhedral set, then $\H^{d-1} ((\p P\cap \O) \smallsetminus
(\p^* P\cap \O)) =0$. We then
obtain that
$$ \widetilde{\phi_{\O}} \,\geq\, \nu_{\min} \times \inf\{ \H^{d-1}(\S \cap \O) \,|\, \S \textrm{ hypersurface that cuts }\G^1
\textrm{ from } \G^2 \textrm{ in }\overline{ \O}\,,\, d(\S, \G^1 \cup \G^2)
>0\}\,.  $$
We recall that the hypersurface $\S$ cuts $\G^1$ from $\G^2$ in $\overline{\O}$ if
$\S$ intersects any continuous path from a point in $\G^1$ to a point in
$\G^2$ that is included in $\overline{\O}$.
We consider such a hypersurface $\S \subset \RR^d$, and we want to bound
from below the quantity $\H^{d-1} (\S \cap \O)$ independently on $\S$.

The idea of the proof is the following. We consider a path from
$\G^1$ to $\G^2$ in $\O$. We construct a tubular
neighbourhood of this path of diameter depending only on the domain
and not on the path itself that lies in $\O$ except at its
endpoints. Then we prove that it is not very deformed compared to a
straight tube. Since $\S$ has to cut this tube, we obtain the desired
lower bound $\H^{d-1} (\S \cap \O)$.

For $i=1,2$, we can find $x_i$ in $\G^i$ and $r_i>0$ such that $\G \cap
B(x_i,r_i) \subset \G^i$ and $\G \cap B(x_i, r_i)$ is a $\C^1$ hypersurface.
We denote by
$v_{\O}(x_i)$ the exterior normal
unit vector to $\O$ at $x_i$, and by $T_{\O} (x_i)$ the hyperplane tangent
to $\G$ at $x_i$. Since $\G$ is of class $\C^1$ in a neighbourhood of 
$x_i$ and $\O$ is a Lipschitz domain, applying lemma \ref{chapitre7rfpres}, we know that for all $\theta>0$, there exists $\eps>0$ depending on
$(\O, \G, \G^1, \G^2, x_1, x_2)$ such that for $i=1,2$ we have 
$$ \left\{ \begin{array}{l} \O \cap B(x_i, 2\eps) \textrm{ is connected}\,,\\
\G \cap B(x_i, 2 \eps) \,\subset \,
    \V_2(T_{\O}(x_i),  2 \eps \sin \theta)  \cap B(x_i,2 \eps) \,,\\
\G \cap B(x_i, 2\eps) \,\subset\, \G^i \,.
\end{array} \right. $$
We fix $\theta$ small enough to have $2\eps \sin \theta < \eps /2$. We
define
$$ A_i \,=\, T_{\O} (x_i) \cap B(x_i, \eps) \quad \textrm{and} \quad D_i
\,=\, \cyl(A_i, \eps) \,, $$
and then
$$ \widehat{\O} = \O \cup \mathring{D}_1 \cup \mathring{D}_2 \,,$$
where $\mathring{D}_i$ is the interior of $D_i$ for $i=1,2$.
We define
$$ X_i\,=\, \{ z \in \mathring{D}_i \,|\, x_iz \cdot v_{\O}(x_i) > \eps/2 \}
\,\subset\, \widehat{\Omega} \,. $$
Then $X_i \subset \widehat{\O} \smallsetminus \O$. Each path $r$
from a point $y_1 \in X_1$ to a point $y_2 \in X_2 $ contains
a path $r'$ from a point $y_1' \in \G^1$ to a point $y_2' \in \G^2$ such
that $r' \subset \overline{\O}$, thus $\S$
intersects $r$. We consider the set
$$ V_i \,=\, \{ z \in X_i \,|\,d_2(z, \partial X_i) > \eps/8  \} \,. $$
Let $\hat{y}_1 \in V_1$, $\hat{y}_2 \in V_2$ such that $d_2(\hat{y}_i, \p X_i)
>\eps / 4$ for $i=1,2$. Since $\widehat{\O}$ is obviously connected by arc,
there exists a path
$\hat{r}$ from $\hat{y}_1$ to $\hat{y}_2$ in $\widehat{\O}$. The path $\hat{r}$ is
compact and $\widehat{\O}$ is open, so $\delta = d_2(\hat{r}, \p \widehat{\O})>0$. We thus can
find a path $r$ included in $\V_2 (\hat{r}, \min (\delta / 2, \eps / 8))$
which is a $\C^\infty$ submanifold of $\RR^d$ of dimension $1$ and which
has one endpoint, denoted by $y_1$, in $V_1$, and the other one, denoted by
$y_2$, in $V_2$.

\begin{figure}
\centering
\begin{picture}(0,0)%
\includegraphics{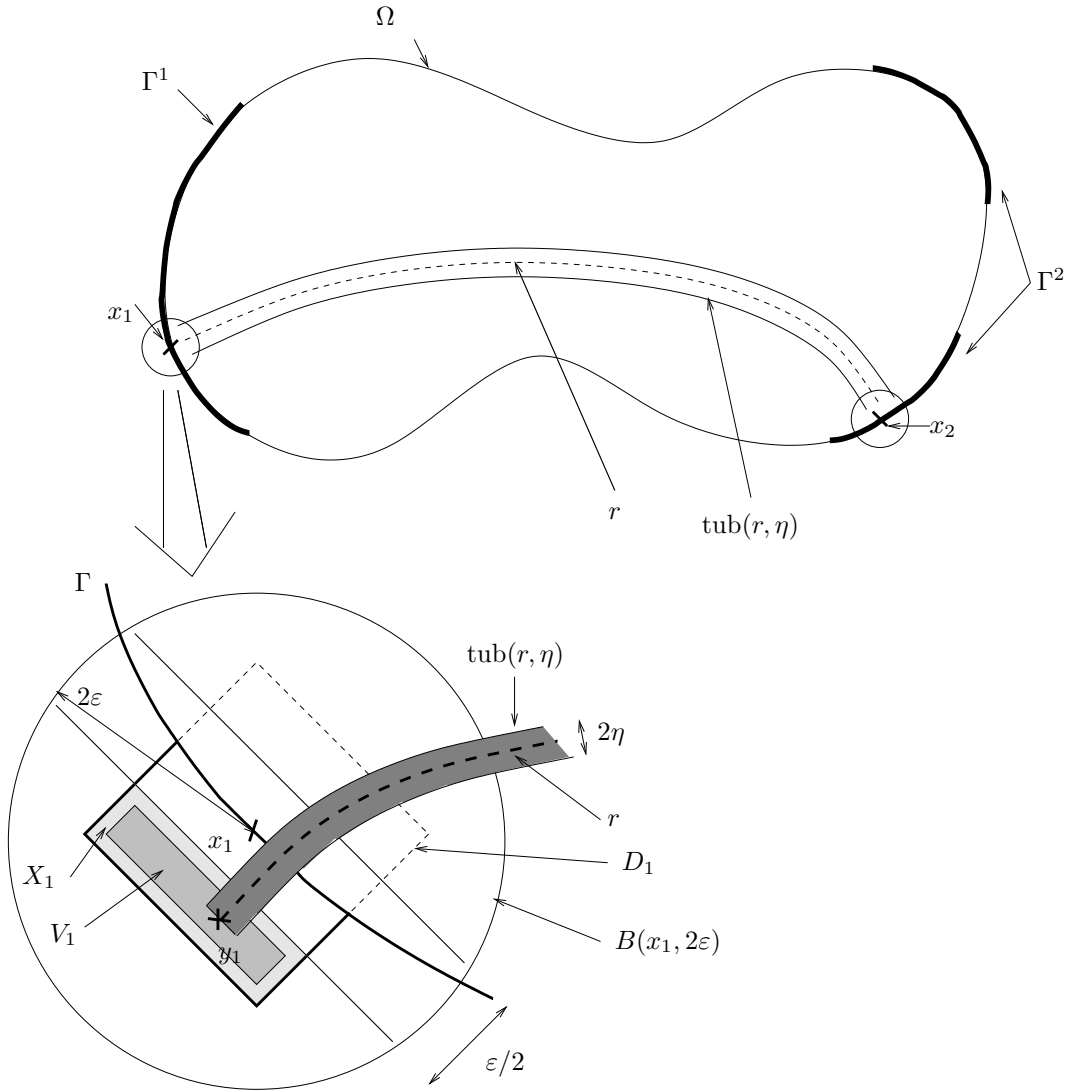}%
\end{picture}%
\setlength{\unitlength}{2368sp}%
\begingroup\makeatletter\ifx\SetFigFont\undefined%
\gdef\SetFigFont#1#2#3#4#5{%
  \reset@font\fontsize{#1}{#2pt}%
  \fontfamily{#3}\fontseries{#4}\fontshape{#5}%
  \selectfont}%
\fi\endgroup%
\begin{picture}(10797,11388)(394,-12643)
\put(1276,-7411){\makebox(0,0)[rb]{\smash{{\SetFigFont{10}{12.0}{\rmdefault}{\mddefault}{\updefault}{\color[rgb]{0,0,0}$\G$}%
}}}}
\put(1126,-11086){\makebox(0,0)[rb]{\smash{{\SetFigFont{10}{12.0}{\rmdefault}{\mddefault}{\updefault}{\color[rgb]{0,0,0}$V_1$}%
}}}}
\put(6826,-10336){\makebox(0,0)[lb]{\smash{{\SetFigFont{10}{12.0}{\rmdefault}{\mddefault}{\updefault}{\color[rgb]{0,0,0}$D_1$}%
}}}}
\put(901,-10486){\makebox(0,0)[rb]{\smash{{\SetFigFont{10}{12.0}{\rmdefault}{\mddefault}{\updefault}{\color[rgb]{0,0,0}$X_1$}%
}}}}
\put(2101,-2161){\makebox(0,0)[rb]{\smash{{\SetFigFont{10}{12.0}{\rmdefault}{\mddefault}{\updefault}{\color[rgb]{0,0,0}$\Gamma^1$}%
}}}}
\put(4351,-1486){\makebox(0,0)[b]{\smash{{\SetFigFont{10}{12.0}{\rmdefault}{\mddefault}{\updefault}{\color[rgb]{0,0,0}$\Omega$}%
}}}}
\put(11176,-4261){\makebox(0,0)[lb]{\smash{{\SetFigFont{10}{12.0}{\rmdefault}{\mddefault}{\updefault}{\color[rgb]{0,0,0}$\Gamma^2$}%
}}}}
\put(5401,-12436){\makebox(0,0)[lb]{\smash{{\SetFigFont{10}{12.0}{\rmdefault}{\mddefault}{\updefault}{\color[rgb]{0,0,0}$\eps /2 $}%
}}}}
\put(1126,-8611){\makebox(0,0)[lb]{\smash{{\SetFigFont{10}{12.0}{\rmdefault}{\mddefault}{\updefault}{\color[rgb]{0,0,0}$2\eps$}%
}}}}
\put(6751,-11161){\makebox(0,0)[lb]{\smash{{\SetFigFont{10}{12.0}{\rmdefault}{\mddefault}{\updefault}{\color[rgb]{0,0,0}$B(x_1,2 \eps)$}%
}}}}
\put(2776,-10111){\makebox(0,0)[rb]{\smash{{\SetFigFont{10}{12.0}{\rmdefault}{\mddefault}{\updefault}{\color[rgb]{0,0,0}$x_1$}%
}}}}
\put(2731,-11266){\makebox(0,0)[b]{\smash{{\SetFigFont{10}{12.0}{\rmdefault}{\mddefault}{\updefault}{\color[rgb]{0,0,0}$y_1$}%
}}}}
\put(5716,-8161){\makebox(0,0)[b]{\smash{{\SetFigFont{10}{12.0}{\rmdefault}{\mddefault}{\updefault}{\color[rgb]{0,0,0}$\tub(r, \eta)$}%
}}}}
\put(8176,-6811){\makebox(0,0)[b]{\smash{{\SetFigFont{10}{12.0}{\rmdefault}{\mddefault}{\updefault}{\color[rgb]{0,0,0}$\tub(r,\eta)$}%
}}}}
\put(6751,-6661){\makebox(0,0)[b]{\smash{{\SetFigFont{10}{12.0}{\rmdefault}{\mddefault}{\updefault}{\color[rgb]{0,0,0}$r$}%
}}}}
\put(1726,-4561){\makebox(0,0)[rb]{\smash{{\SetFigFont{10}{12.0}{\rmdefault}{\mddefault}{\updefault}{\color[rgb]{0,0,0}$x_1$}%
}}}}
\put(10051,-5761){\makebox(0,0)[lb]{\smash{{\SetFigFont{10}{12.0}{\rmdefault}{\mddefault}{\updefault}{\color[rgb]{0,0,0}$x_2$}%
}}}}
\put(6676,-9871){\makebox(0,0)[lb]{\smash{{\SetFigFont{10}{12.0}{\rmdefault}{\mddefault}{\updefault}{\color[rgb]{0,0,0}$r$}%
}}}}
\put(6571,-8971){\makebox(0,0)[lb]{\smash{{\SetFigFont{10}{12.0}{\rmdefault}{\mddefault}{\updefault}{\color[rgb]{0,0,0}$2\eta$}%
}}}}
\end{picture}%
\caption{Construction of $\tub(r,\eta)$.}
\label{chapitre7positivite}
\end{figure}

As we explained previously, $d_2(r, \p \widehat{\O})>0$, so there exists a positive
$\eta_1$ such that $\V_2(r, \eta_1) \subset \widehat{\O}$. We can suppose that $\eta_1 <
\eps/16$, to obtain that $B(y_i, \eta_1) \subset X_i$ for $i=1,2$. For all $z$ in $r$ we denote by $N_r(z)$ the hyperplane orthogonal to
$r$ at $z$, and by $N^{\eta}_r (z)$ the subset of $N_r(z)$ composed of the
points of $N_r(z)$ that are at distance smaller than or equal to $\eta$ of $z$.
The tubular neighbourhood of $r$ of radius $\eta$, denoted by
$\tub(r,\eta)$, is the set of all the points $z$ in $\RR^d$ such that there
exists a geodesic of length smaller than or equal to $\eta$ from $z$ that
meets $r$ orthogonally, i.e.,
$$ \tub(r, \eta) \,=\, \bigcup_{z\in r} N_r^{\eta} (z)\,, $$
(see for example \cite{Gray}). We have a picture of this tubular
neighbourhood on figure \ref{chapitre7positivite}.
Since $r$ is a compact $\C^\infty$ submanifold of $\RR^d$ which is complete, there exists a $\eta_2>0$
small enough such that for all $\eta \leq \eta_2$, the tubular neighbourhood
of $r$ of diameter $\eta$ is well defined by a $\C^\infty$-diffeomorphism (see
for example \cite{BergerGostiaux}, Theorem 2.7.12, or \cite{Gray}), i.e.,
there exists a $\C^{\infty}$-diffeomorphism $\psi$ from 
$$Nr^{\eta}\,=\,\{(z,v)\,,\, z\in r\,,\, v\in N^{\eta}_r(z)\} $$
to $\tub(r,\eta)$. We choose a positive $\eta$ smaller than $\min (\eta_1,
\eta_2)$. We stress the fact that this $\eta$ depends on $(\O, \G, \G^1,
\G^2)$ but not on $\S$.

Let $(I,h)$ be a parametrisation
of class $\C^{\infty}$ of $r$, i.e., $I=[a,b]$ is a closed interval of
$\RR$, $h: I\rightarrow r $ is a $\C^{\infty}$-diffeomorphism
which is an immersion. Let $z$ be in $r$, and $u_z = h^{-1}(z) \in I$. The
vector $h'(u_z)$ is tangent to $r$ at $z$, and there exists some vectors
$(e_2(z),...,e_d(z))$ such that $(h'(u_z),e_2(z),...,e_d(z))$ is a direct
basis of
$\RR^d$. There exists a neighbourhood $U_z$ of $u_z$ in $I$ such that for
all $u\in U_z$, $(h'(u),e_2(z),...,e_d(z))$ is still a basis of $\RR^d$,
since $h'$ is continuous. Indeed the condition for a family of vectors
$(\alpha_1,...,\alpha_d)$ to be a basis of $\RR^d$ is an open condition,
because it corresponds to $\deter ((\alpha_1,...,\alpha_d))>0$ where $\deter$
is the determinant of the matrix. We apply the Gram-Schmidt process to the basis $(h'(u),e_2(z),...,e_d(z))$ to obtain a direct orthonormal basis
$(h'(u)/\| h'(u) \|, v_2(u,z),...,v_d(u,z))$ of $\RR^d$ for all $u \in
U_z$, such that the dependence of $(h'(u)/\| h'(u) \|, v_2(u,z),...,v_d(u,z))$
on $u\in U_z$ is of class $\C^\infty$. We remark that the family
$(v_2(u,z), ..., v_d(u,z))$ is a direct orthonormal basis of $N_r(h(u))$ for all
$u\in U_z$. We have associated with each $z\in r$ a neighbourhood $U_z$ of
$u_z = h^{-1}(z)$ in $I$, we can obviously suppose that $U_z$ is an
interval which is open in $I$. Since $(U_z, z\in r)$ is a covering of the
compact $I$, we can extract a finite covering $(U_j, j= 1,...,n )$ from
it. We can choose this family to be minimal, i.e., such that $(U_j, j\in
\{1,...,n\}
\smallsetminus j_0)$ is not a covering of $I$ for any $j_0 \in
\{1,...,n\}$. We then reorder the $(U_j, j =1,...,n)$ (keeping the same
notation) by the increasing order of their left end point in $I\subset
\RR$. Since the family $(U_j)$ is minimal, each point of $I$ belongs either
to a unique set $U_j$, $j\in \{1,...,n\}$, or to exactly two sets $U_j$ and
$U_{j+1}$ for $j\in \{1,...,n-1 \}$.
We denote by $a_j$ the middle of the non-empty open
interval $U_j \cap U_{j+1}$ for $j\in \{1,...,n-1 \}$, and by $(h'(u)/\| h'(u) \|,v_2(u,j),...,v_d(u,j))$ the direct orthonormal basis defined previously on
$U_j$ for $j\in \{1,...,n \}$. We want to construct a family of direct
orthonormal basis $(h'(u)/\|h'(u)\|, f_2(u),...,f_d(u))$ of $\RR^d$ such
that the function:
$$\psi: u \in I \mapsto (h'(u)/\|h'(u)\|, f_2(u),...,f_d(u)) $$
is of class $\C^\infty$. We have to define a
concatenation of the $(h'(u)/\| h'(u) \|, v_2(u,j),...,v_d(u,j))$ over the
different sets $U_j$. For $u\in [a, a_1] $, we
define
$$\psi(u) \,=\, (h'(u)/\| h'(u) \|,
v_2(u,1),...,v_d(u,1)) \,.$$
Thus the function $\psi$ defined on $[a, a_1]$ is of class
$\C^{\infty}$. On $U_1 \cap U_2$ we have defined two different
direct orthonormal basis $(h'(u)/\| h'(u) \|,
v_2(u,j),...,v_d(u,j))$ for $j=1$ and $j=2$
that have the same first vector. Let $\phi_1 : U_1 \cap U_2 \rightarrow
SO_{d-1} (\RR) $ be the function of class $\C^{\infty}$ that associates to
each $u\in U_1\cap U_2$ the matrix of change of basis from $
(v_2(u,2),...,v_d(u,2))$ to $ (v_2(u,1),...,v_d(u,1))$.


If $b_1$ is the right end point of $U_1 \cap U_2$, then $\phi_1$ is in
particular defined on $[a_1,b_1[$. Let $g_1$ be a $\C^\infty$-diffeomorphism
from $[a_1, b_1[$ to $[a_1, \infty[$ which is strictly increasing (so
$g_1(a_1) = a_1$) and such that all the derivatives of $g_1$
at $a_1$ are null. Then $\phi_1 \circ g_1^{-1}$ is defined on $[a_1,
+\infty[$ and all its derivatives at $a_1$ are equal to those of $\phi_1$. We
then transform all the orthonormal basis $( v_2(u,j),...,v_d(u,j))$ of
$\RR^{d-1}$ for $j\geq 2$ and $u\geq a_1$ by the change of basis $\phi_1
\circ g_1^{-1}$, and we denote the new direct orthonormal basis of $\RR^{d-1}$ obtained
this way by $( \widetilde{v}_2(u,j),...,\widetilde{v}_d(u,j))$. We then
define $\psi$ on $]a_1, a_2]$ by
$$ \psi(u) \,=\,
(h'(u)/\|h'(u)\|,\widetilde{v}_2(u,2),...,\widetilde{v}_d(u,2)  )  \,,$$
and we remark that $\psi(u)$ still defines a direct orthonormal basis of $\RR^d$.
The function $\psi$ is of class $\C^{\infty}$ on $[a, a_2]$, including at
$a_1$. We iterate this process with the family of basis
$$(h'(u)/\|h'(u)\|,\widetilde{v}_2(u,j),...,\widetilde{v}_d(u,j)  )
  \,,\,j= 2,..., n $$
at $a_2$, etc..., finitely many times since we work with a
finite covering of $I$. We obtain in the end a function 
$$ \psi \circ h^{-1} : r \rightarrow SO_{d-1}(\RR)$$
which is of class $\C^\infty$, and for all $z\in r$, the set of the points of
$\RR^d$ that have for first coordinate $0$ in the basis $\psi \circ h^{-1}
(z)$ is exactly the hyperplane $N_r(z)$.

For each $t=(t_2,...,t_{d-1}) \in \{ z\in \RR^{d-1} \,|\, d(z, 0) \leq \eta  \}$, the set
$$ r_t \,=\, \{ y\in \RR^d \,|\, \exists z\in r\,,\,y \textrm{
  has coordinates } (0, t_2,...,t_{d-1}) \textrm{ in the basis } \psi\circ
h^{-1} (z)\} $$
is a continuous path (even of class $\C^\infty$) from a point in $X_1$ to a
point in $X_2$, therefore
$$ r_t \cap \S \cap \overline{\O} \,\neq \,\varnothing \,.  $$
Moreover, since $d(\S, \G^1\cup \G^2)>0$, we obtain that
\begin{equation}
\label{chapitre7truc}
r_t \cap \S \cap \O \,\neq \,\varnothing \,.
\end{equation} 
For each $y\in \tub
(r, \eta)$, there exists a unique $z_y\in r$ such that $y\in N_r(z_y)$, so we
can associate to $y$ its coordinates $(0,t_2(y),...,t_d(y))$ in the
basis $\psi \circ h^{-1} (z_y)$. We define the
projection $p$ of $\tub (r,\eta)$ on $N_r^{\eta}(y_1)$ that associates to each
$y$ in $\tub(r,\eta)$ the point of coordinate $(0, t_2(y),...,t_d(y))$ in
the basis $\psi \circ h^{-1} (y_1)$. Then $p$ is of class $\C^\infty$ as is
$\psi \circ h^{-1}$. If $z$ belongs to $ N_r^\eta (y_1)$, and $t(z) =
(t_2(z),...,t_d(z)) $, then we know by equation (\ref{chapitre7truc}) that
there exists a point on $r_{t(z)}$ that intersects $\S$ in $\O$. Moreover, $r_{t(z)}$ is
exactly the set of the points $y$ of $\tub (r,\eta)$ whose image $p(y)$ by
this projection is the point $z$. Thus
$$ p(\S \cap \tub(r, \eta) \cap \O) \,=\, N_r^\eta (y_1)\,. $$
Since $\tub (r, \eta)$ is compact, $p$ is a Lipschitz
function on $\tub(r, \eta)$, and so there exists a constant $K$, depending on
$p$, hence on $\O$, $r$, $\eta$, but not on $\S$, such that
$$  \H^{d-1} (\S \cap \O) \,\geq\, \H^{d-1} (\S \cap \tub (r, \eta ) \cap \O) \,\geq\, K
\H^{d-1} (p(\S \cap \tub (r, \eta))) \,\geq\, K \alpha_{d-1} \eta^{d-1} \,. $$
This ends the proof of the positivity of $\widetilde{\phi_{\O}}$ when
$\Lambda (0) < 1-p_c(d)$.\\

\noindent
{\bf Acknowledgment:} The second author would like to warmly thank Immanuel
Halupczok, Thierry L\'evy and Fr\'ed\'eric Paulin for helpful discussions.


\def\cprime{$'$}

\end{document}